\documentclass[reqno, 11pt,twoside, final]{amsart}
\pagespan{1}{\pageref*{LastPage}}
\usepackage{etoolbox,lastpage}
\commby{}
\usepackage{amsmath,amsthm,amscd,amsfonts,amssymb,enumerate}
\usepackage{graphicx}
\usepackage{color}
\usepackage[colorlinks]{hyperref}
\numberwithin{equation}{section}

\begin{document}

\title[ ]{Subdifferential of fuzzy n-cell number valued functions and its applications in optimization problems}

\author[S. Fatemi, I. Sadeqi and F. Moradlou]{Samira Fatemi$^1$, Ildar Sadeqi$^{*2}$ and Fridoun Moradlou$^3$}
\address{\indent $^{1,2,3}$ Department of Mathematics
\newline \indent Sahand University of Technology
\newline \indent Tabriz, Iran}
\email{\rm $^1$ s\_fatemi98@sut.ac.ir \& samira.fatemi09@gmail.com}
\email{\rm $^2$ esadeqi@sut.ac.ir \& esadeqi@yahoo.com}
\email{\rm $^2$ moradlou@sut.ac.ir \& fridoun.moradlou@gmail.com}

 \thanks{$^*$Corresponding author}

%
 \maketitle
 \newtheorem{theorem}{Theorem}[section]
\newtheorem{proposition}[theorem]{Proposition}
\newtheorem{lemma}[theorem]{Lemma}
\newtheorem{corollary}[theorem]{Corollary}
\theoremstyle{definition}
\newtheorem{definition}[theorem]{Definition}
\newtheorem{example}[theorem]{Example}
\newtheorem{remark}[theorem]{Remark}
\newtheorem{fact}[theorem]{Fact}
%
\begin{abstract}
In this paper, we present the concept of subdifferential for fuzzy n-cell number valued functions. Then we state some theorems related to subdifferentiability based on the new definition. Finally, we present some applications emphasized on optimization problems, including the Lagrangian dual problem and minimizing the composite problem.


\end{abstract}
\section{\bf Introduction}
The theory of fuzzy sets was introduced by Zadeh in 1965, and so far, many mathematicians have worked on the theory of fuzzy analysis and its applications, which improved this field of mathematics.
The importance of studying fuzzy analysis is clearly known from both theoretical and practical points of view. As known the definition of a fuzzy subset is obtained by extending the range from $\{0,1\}$ to the whole of $[0,1].$ Therefore, the classical set is a special case of a fuzzy set. In fact, fuzzy numbers are a powerful tool for modeling uncertainty and processing vague or subjective information in mathematical models.
The space of fuzzy numbers is not a vector space because it does not contain the additive inverse of its members. That is, for every fuzzy number $u,$ $u-u \neq \{0\},$ unless $u$ is single-membered.
This defect makes it difficult to study mathematical concepts in many cases. For example, fuzzy differential equations. To solve this problem, Hukuhara introduced a difference between two fuzzy sets so that based on it, the space of fuzzy numbers will include the reciprocal of each member \cite{Hukuhara}.
The Hukuhara difference does not always exist, so Stefanin introduced the generalized Hukuhara difference \cite{Stefanini1,Stefanini2}. The latter difference may not always exist. Therefore, Bede and Stefanin introduced a new generalization of Hukuhar's difference \cite{Bede}. Based on Hukuhara difference and its generalizations, derivative of fuzzy functions \cite{Bede,Hai,Hai1,Qiu,Gong}, directional derivative and subdifferential of fuzzy functions were defined \cite{Gong1,Wang1,Gong}. The set of fuzzy numbers is a partially ordered set. Therefore, two fuzzy numbers may not be comparable. Based on existing partial order relations and different types of fuzzy numbers, many researches provided ways to solve optimization problems using derivative and subdifferential of fuzzy functions \cite{Gong1,Qiu,Wang1,Gong}.
Since the subdifferential is a powerful tool for solving fuzzy optimization problems, our purpose in this paper is to provide a definition of the subdifferential of fuzzy functions and apply it to the dual Lagrange problems and Karush–Kuhn–Tucker optimality conditions. Accordingly, in the third section of this paper, first the definitions of derivative, directional derivative, gradient of n-cell fuzzy functions and related theorems are discussed. Then, the definition of subdifferential of n-cell fuzzy functions is given. Finally, fuzzy optimization problems and their applications are studied.
\section{\bf Preliminaries}
A fuzzy subset $u$ is a function as $u:\mathbb{R}^n \to [0,1].$ For each fuzzy set $u,$ we denote the $r$-level set of $u$ by $ [u]^r=\{ x \in \mathbb{R}^n: u(x) \geq r \}$ for any $ r \in (0,1]$ and the support of $u$ by $ \text{supp} \ u=\{ x \in \mathbb{R}^n: u(x) > 0 \}.$
\begin{definition}\cite{Wu}
Let $ u:\mathbb{R}^n \to [0,1] $ satisfies the following properties:
\begin{itemize}
\item[(1)]
$u$ is normal; that is, there is an $x_0\in \mathbb{R}^n $ with $ u(x_0)=1$;
\item[(2)]
$u$ is a convex fuzzy set; that is, $ u(\lambda x +(1-\lambda)y) \geqslant \min (u(x),u(y)) $ whenever $x,y\in \mathbb{R}^n$ and $ \lambda \in [0,1]$;
\item[(3)]
$u(x)$ is upper semi-continuous, i.e., $ u(x_0) \geq \overline{\lim}_{k \to \infty} u(x_k)$ for any $x_k \in \mathbb{R}^n \ (k=0,1,2,...), \ x_k \to x_0;$
\item[(4)]
$ \overline{\text{supp} \ u}=\overline{ \{x\in \mathbb{R}^n: u(x)>0 \}}$ is a compact set.\\
\end{itemize}
Then, $u$ is called the fuzzy number. The fuzzy numbers space is denoted by $E^n.$
\end{definition}
Each $a \in \mathbb{R}^n$ is a fuzzy number defined by
\begin{equation*}
 \hat{a} (x) = \left\{
\begin{array}{cc}
1 & \text{if } x=a,\\
0 & \text{if } x \neq a,
\end{array} \right.
\end{equation*}
for $ x \in \mathbb{R}^n.$ The addition and scalar multiplication operations on $E^n$ are defined as follow:
\begin{align*}
 (u+v)(x) =\sup_{y+z=x} \min \{u(y), v(z) \},\\
\end{align*}
 \begin{align*}
(\lambda u)(x)= \left\{
\begin{array}{cc}
u(\lambda^{-1} x)  & \text{if } \lambda \neq 0,\\
\hat{0} & \text{if } \lambda = 0,
\end{array} \right.
 \end{align*}
for $ u,v \in E^n$ and $\lambda \in \mathbb{R}.$
\begin{definition}\cite{Wu}
Let $ u \in E^n$ such that
$$ [u]^r=\prod_{i=1}^n[u^-_i(r),u^+_i(r)]$$
for $r\in[0,1],$ where $ u^-_i(r),u^+_i(r) \in \mathbb{R}$ with $u^-_i(r)\leq u^+_i(r), \ i=1,2,...,n,$ then $u$ is said to be a fuzzy n-cell number. The set of all fuzzy $n$-cell numbers is denoted by $L(E^n).$
\end{definition}
\begin{theorem} \cite{Wu}  \label{l}
Let $u,v\in L(E^n)$ and $k\in \mathbb{R}.$ Then, for any $r\in [0,1],$\\\\
 $ (1) [u+v]^r=[u]^r+[v]^r=\prod_{i=1}^n[u^-_i(r)+v^-_i(r),u^+_i(r)+v^+_i(r)],$\\
 \begin{equation*}
\hspace*{-6.7cm} (2) [ku]^r=k[u]^r=\left\{
\begin{array}{rl}
\prod_{i=1}^n [ku^-_i(r),ku^+_i(r)] & k\geq 0,\\
\prod_{i=1}^n [ku^+_i(r),ku^-_i(r)] & k< 0,
\end{array} \right.
\end{equation*}
 $ (3) [u.v]^r=\prod_{i=1}^n \big[\min \{u^-_i(r) v^-_i(r),u^-_i(r) v^+_i(r),u^+_i(r) v^-_i(r),u^+_i(r) v^+_i(r)\},$\\
$ \hspace*{3cm} \max \{u^-_i(r) v^-_i(r),u^-_i(r) v^+_i(r),u^+_i(r) v^-_i(r),u^+_i(r) v^+_i(r)\} \big].$
\end{theorem}
\begin{remark}\cite{Wu}
We denote
$$LC(\mathbb{R}^n)=\{A| \ \text{there exists} \ a_i \leq b_i  \ \text{such that} \ A=\prod_{i=1}^n [a_i,b_i], i=1,2,...,n \}.$$
\end{remark}
\begin{definition}\label{35} \cite{Wu}
A mapping $d_L:LC(\mathbb{R}^n) \times LC(\mathbb{R}^n) \to [0,+\infty)$ is defined as
\begin{align*}
d_L(A,B) &=\max_{1\leq i \leq n}  \Big \{ |a^-_i-b^-_i|,|a^+_i -b^+_i| \Big \}\\
&=\max \Big \{ |a^-_1-b^-_1|,|a^+_1 -b^+_1|,...,|a^-_n-b^-_n|,|a^+_n -b^+_n| \Big \}
\end{align*}
for any $A=\prod_{i=1}^n [a^-_i,a^+_i],$ $B=\prod_{i=1}^n [b^-_i,b^+_i].$
\end{definition}
\begin{theorem} \cite{Wu} \label{z}
For any $ A,B,C \in LC(\mathbb{R}^n),$ and each $\alpha \in \mathbb{R},$ $d_L$ satisfies
\begin{itemize}
\item[(1)] $d_L(A,B)=d_L(B,A),$
\item[(2)] $d_L(A,B)\geq 0,$
\item[(3)] $d_L(A,B)=0 \Leftrightarrow A=B,$
\item[(4)] $d_L(A,B) \leq d_L(A,C)+d_L(C,B),$
\item[(5)] $d_L(A+C,B+C)=d_L(A,B),$
\item[(6)] $d_L(\alpha A,\alpha B)=|\alpha|d_L(A,B),$
\end{itemize}
where $A=\prod_{i=1}^n [a^-_i,a^+_i],$ $B=\prod_{i=1}^n [b^-_i,b^+_i]$ and $C=\prod_{i=1}^n [c^-_i,c^+_i].$
\end{theorem}
\begin{theorem} \cite{Wu} \label{cc}
The mapping $D_L:L(E^n) \times L(E^n) \to [0,+\infty)$ is defined by $D_L(u,v)=\sup_{r \in [0,1]} d_L([u]^r,[v]^r)$ for any $(u,v) \in L(E^n) \times L(E^n).$ For each $u,v,w \in L(E^n)$ and any $\alpha \in \mathbb{R},$ $D_L$ satisfies
\begin{itemize}
\item[(1)] $D_L(u,v)=D_L(v,u),$
\item[(2)] $D_L(u,v)\geq 0,$
\item[(3)] $D_L(u,v)=0 \Leftrightarrow u=v,$
\item[(4)] $D_L(u,v) \leq D_L(u,w)+D_L(w,v),$
\item[(5)] $D_L(u+w,v+w)=D_L(u,v),$
\item[(6)] $D_L(\alpha u,\alpha v)=|\alpha|D_L(u,v).$
\end{itemize}
\end{theorem}
\begin{theorem}  \cite{Wu}
$(L(E^n),D_L)$ is a complete metric space.
\end{theorem}
From now on, the set $M$ is supposed to be a nonempty subset of $\mathbb{R}^m.$
\begin{definition} \cite{Hai}\label{18}
Suppose that $ F:M \to L(E^n)$ is a fuzzy $n$-cell number-valued function. If for any $ \epsilon>0$, exists $\delta>0$ such that for all $t \in M, ||t-t_0||<\delta$ ensures that
$$ D_L(F(t),u)<\epsilon,$$
then $\lim_{t \to t_0} F(t)=u.$ The function of $F$ is said to be continuous at $t_0$ if and only if $ \lim_{t\to t_0} F(t) =F(t_0). $
\end{definition}
\begin{definition}\label{aa}
Let $\{ u^{(k)} \}$ be a fuzzy $n$-cell sequence, then $ \lim_{k \to \infty} { u^{(k)} }=u,$ if and only if for any $\epsilon> 0$ exists an integer $ N>0,$ such that
$$ D(u^{(k)},u)<\epsilon,$$
 whenever $ k \geq N.$
\end{definition}
\begin{definition} \cite{Wang}\label{c}
For any  $u,v \in L(E^n),$ define $u\leq v$ if and only if $u^-_i (r)\leq v^-_i (r)$ and $u^+_i (r) \leq v^+_i (r), \ i = 1, 2, . . . , n,$ for any $r \in [0,1]$ (denoted by $[u]^r \leq [v]^r $ for any $r \in [0,1]).$
\end{definition}
A subset $W \in  L(E^n)$ is said to be order bounded from above if there exists a fuzzy $n$-cell number $s \in L(E^n),$ called an upper bound of $W,$ such that $u\leq s$ for every $u \in W.$ A fuzzy $n$-cell number $s$ is called the supremum of $W$ if $s$ is an order upper bound of $W$ and $s\leq s'$ holds for any upper bound $s'$ of $W,$ so we denote $ s=\sup_{u\in W} u$ (or $ s=\sup W$). An order lower bound and the infimum of $W$ can be defined similarly. The set $W$ is said to be order bounded if it is both order bounded from above and order bounded from below.
\begin{theorem}\label{r}\cite{Wu,Wang}
If $ u \in L(E^n),$ then
\begin{itemize}
\item[(1)] $[u]^r$ is a non-empty n-dimensional closed polyhedron, i.e., $ [u]^r=\prod_{i=1}^n [u_i^-(r),u_i^+ (r)]$ for each $ r \in [0,1];$
\item[(2)] $ [u]^{r_2} \subset [u]^{r_1}$ for $ 0\leq r_1 \leq r_2 \leq 1;$
\item[(3)] if $ \{r_n\} $ is a positive non-decreasing sequence converging $ r \in (0,1],$ then $\bigcap_{n=1}^\infty [u]^{r_n}=[u]^r.$
\end{itemize}
Conversely, if for each $ r\in [0,1],$ $A_r$ is an element of $K(\mathbb{R}^n)$ which satisfies Conditions $(1)-(3),$ then
there exists a unique $n$-dimensional fuzzy number $u \in L(E^n)$ such that $[u]^r=A_r$ for each $ r\in(0,1]$ and $[u]^0=\overline{\cup_{r \in (0,1]} A_r} \subset A_0.$
\end{theorem}
\begin{theorem}\label{f}\cite{Wu,Wang}
If $u \in L(E^n),$ then for $i = 1, 2, . . . , n,$ $u^-_i (r), u^+_i(r)$ are real-valued functions on $[0,1],$ and satisfy
\begin{itemize}
\item[(1)] $u^-_i(r)$ are non-decreasing and left continuous;
\item[(2)] $u^+_i(r)$ are non-increasing and left continuous;
\item[(3)] $u^-_i(r)\leq u^+_i(r)$ (it is equivalent to $u^-_i(1)\leq  u^+_i(1));$
\item[(4)] $u^-_i(r), u^+_i(r)$ are right continuous at $r=0.$
\end{itemize}
Conversely, if $a_i(r), b_i(r), i = 1, 2, . . . , n,$ are real-valued functions on $[0,1]$ which satisfy Conditions $(1)-(4),$ then there exists a unique $u \in L(E^n)$ such that $ [u]^r=\prod_{i=1}^n [a_i(r),b_i(r)]$ for any $ r \in [0,1].$
\end{theorem}
\begin{definition} \cite{Bede}
The generalized Hukuhara difference ($gH$-difference for short) $u \ominus_g v$ is defined as:
$$ u \ominus_{gH} v=w \Leftrightarrow u=w+v \quad or \quad v=u+(-1)w.$$
for any $u,v \in E^n.$
\end{definition}
\begin{definition} \cite{Gomes1}
The generalized difference ($g$-difference for short) $u \ominus_g v$ is defined by its level sets:
$$[u \ominus_g v]^r=cl \Big( conv \bigcup_{\beta \geq r} ([u]^\beta \ominus_{gH} [v]^\beta) \Big) \ \text{for any $r\in[0, 1]$},$$ where the $gH$-difference is with interval operands $ [u]^\beta $ and $[v]^\beta.$
\end{definition}
\begin{theorem}\label{k} \cite{Hai}
Let $u, v \in L(E^n).$ If the $g$-difference $u\ominus_g v$ of $u$ and $v$  exists, then for any $r\in[0, 1],$ we have
\begin{align*}
[u\ominus_g v]^r &=\prod_{i=1}^n [\inf_{\beta \geq r}\min \{u^-_i(\beta)-v^-_i(\beta),u^+_i(\beta)-v^+_i(\beta)\},\\
 & \hspace*{1.2cm} \sup_{\beta \geq r} \max \{u^-_i(\beta)-v^-_i(\beta),u^+_i(\beta)-v^+_i(\beta)\}].
\end{align*}
\end{theorem}
\begin{remark} \cite{Hai}
A sufficient condition for the existence of $u\ominus_g v$ is that either $[u^-_i(r),u^+_i(r)],$ contains a translation of $ [v^-_i(r),v^+_i(r)]$ or $[v^-_i(r),v^+_i(r)]$ contains a translation of $ [u^-_i(r),u^+_i(r)]$ for any $r\in[0, 1],$ $i=1, 2, ..., n.$ That is, if $l_i[u]^r\geq l_i[v]^r$ or $l_i[u]^r< l_i[v]^r$ for any $r\in[0, 1],$ $i=1, 2, ..., n,$ then the $g$-difference $u\ominus_g v$ exists and $u\ominus_g v \in L(E^n).$
\end{remark}
For any $r\in[0,1],$ $l_i[u]^r=u^+_i(r)-u^-_i(r), \ i=1, 2, ..., n,$ is called the $r$-level length of $u\in L(E^n)$ concerning the $i$th component.
\begin{theorem}\label{bb} \cite{Hai}
If $ u,v,w \in L(E^n),$ and $u \ominus_g v,$ $ v\ominus_g u,$ $u \ominus_g (-v)$ and $(-u) \ominus_g v$ exist, then we have
\begin{itemize}
\item[(1)] $ u \ominus_g u=\hat{0},$ $ u \ominus_g \hat{0}=u,$ $ \hat{0} \ominus_g u=-u,$
\item[(2)] $u \ominus_g v=-(v \ominus_g u),$
\item[(3)] $ k(u\ominus_g v)=ku \ominus_g kv$ for any $k \in \mathbb{R},$
\item[(4)] $k_1 u \ominus_g k_2 u=(k_1-k_2)u$ for any $k_1, k_2 \in \mathbb{R}$ and $ k_1.k_2\geq 0,$
\item[(5)] $ u \ominus_g (-v)=v \ominus_g (-u)$ and $ (-u)\ominus_g v=(-v) \ominus_g u,$
\item[(6)] $ (u+v)\ominus_g v=u,$
\item[(7)] $ \hat{0} \ominus_g (u\ominus_g v)= v \ominus_g u=(-u)\ominus_g (-v),$
\item[(8)] $ u \ominus_g v= v \ominus_g u=w$ if and only if $ w=-w.$
\end{itemize}
\end{theorem}
\begin{theorem}\label{xx} \cite{Wu}
Let $ u,v,w \in L(E^n)$ and $ k, k_1,k_2 \in \mathbb{R}.$ Then $u+v, \ ku, \ uv \in L(E^n)$ and
\begin{itemize}
\item[(1)] $u \leq v$ if and only if $ u+w \leq v+w;$
\item[(2)] If $ u \leq v$ then $ k_1 u\leq k_1 v,$ $k_2 u \geq k_2 v$ for any $k_1 \geq 0,$ $ k_2 < 0;$
\item[(3)] $k(u+v)=ku+kv;$
\item[(4)] $k_1(k_2 u)= (k_1 k_2)u;$
\item[(5)] $ (k_1+k_2) u= k_1 u+k_2 u$ when $ k_1 \geq 0$ and $k_2 \geq 0.$
\end{itemize}
\end{theorem}
\begin{theorem}\label{mm} \cite{Hai}
Let $a,b,c,d \in L(E^n).$ If $l_i[a]^r \geq l_i[c]^r$ and $l_i[b]^r \geq l_i[d]^r$ or $l_i[a]^r < l_i[c]^r$ and $l_i[b]^r < l_i[d]^r$ for any $r \in [0,1], \ i=1,2,...,n,$ then we have
$$ (a+b)\ominus_g (c+d)=(a \ominus_g c)+(b \ominus_g d).$$
\end{theorem}
\begin{definition}\cite{Liu}
Let $u=(u_1,u_2,...,u_m) \in (E^n)^m$ and $t=(t_1,t_2,...,t_m) \in \mathbb{R}^m$ be an n-dimensional fuzzy vector and an $n$-dimensional real vector, respectively. The product of a fuzzy vector with a real vector is defined as $ t.u= \sum_{i=1}^m t_i u_i,$ which is an $n$-dimensional fuzzy number.
\end{definition}
\begin{theorem}\label{12}
Let $F:M \to (-\infty,+\infty)$ be differentiable at $t \in \text{int(domf)}.$ Then  for every $ d \in M$
$$F'(t;d)=\langle\nabla F(t),d\rangle.$$
\end{theorem}
\section{\bf Main Results}
\subsection{\bf Directional derivative of fuzzy functions}
Here, we give the notion of convexity for subsets of $L(E^n)$ and fuzzy $n$-cell number-valued functions. Also, we generalize the concepts of directional derivative for fuzzy $n$-cell number-valued functions.
\begin{definition}
A subset $C$ of $L(E^n)$ is convex, if and only if for any $ u,v \in L(E^n)$ and $ \lambda \in [0,1],$
$$ \lambda u+ (1-\lambda) v \in C.$$
\end{definition}
\begin{definition}
Let $F: M \to L(E^n)$ be a fuzzy function and $M \subseteq \mathbb{R}^m$ be a convex set. The function $F$ is convex on $M$ if and only if
$$ F(\lambda x+(1-\lambda) y)\leq \lambda F(x)+(1-\lambda) F(y),$$
for any $x,y \in M$ and $ \lambda \in [0,1].$
\end{definition}
\begin{theorem}\label{52}
Assume that $F: M \to L(E^n)$ is a fuzzy function, and $M \subseteq \mathbb{R}^m$ is a convex set. Then $F$ is convex on $M$ if and only if $F^-_i(r,t)$ and $F^+_i(r,t),$ $i=1,2,...,n,$ are convex real-valued functions with respect to $t$ for any $r \in [0,1].$
\end{theorem}
\begin{proof}
$F$ is convex if and only if
$$ F(\lambda x+(1-\lambda) y)\leq \lambda F(x)+(1-\lambda) F(y),$$
if and only if according to Definition \ref{c},
$$ F^-_i(r,\lambda x+(1-\lambda) y)\leq \lambda F^-_i(r,x)+(1-\lambda)F^-_i(r,y),$$
$$ F^+_i(r,\lambda x+(1-\lambda) y)\leq \lambda F^+_i(r,x)+(1-\lambda)F^+_i(r,y),$$
for any $ r,\lambda\in [0,1]$ and $ x,y \in M.$
\end{proof}
The following definition is a generalization of Definition 1 of \cite{Chalco}.
\begin{definition}
Let $F: M \to L(E^n)$ be a fuzzy $n$-cell number-valued function and $ t_0 \in M.$
We say that $[L]^r$ is a limit of $[F(t)]^r$ at $t_0,$ if for any $\epsilon > 0,$ exists $\delta >0$ such that for any $t \in M,$ $ ||t-t_0||<\delta$ implies $d_L ([F(t)]^r,[L]^r)< \epsilon,$ for any $r \in [0,1].$ We denote it as
$$ \lim_{t \to t_0} [F(t)]^r=[L]^r,$$
where
$$[F(t)]^r=\prod_{i=1}^n [F_i^-(r,t),F_i^+(r,t)], \quad [L]^r=\prod_{i=1}^n [L_i^-(r),L_i^+(r)].$$
\end{definition}
\begin{theorem}  \label{vv}
Let $ F:M \to L(E^n)$ be a fuzzy $n$-cell number-valued function and $t_0 \in M.$ Then $\lim_{t \to t_0} F_i^-(r,t)$ and $\lim_{t \to t_0} F_i^+(r,t)$ exist if and only if $\lim_{t \to t_0} [F(t)]^r$ exists for any $ i=1,2,...,n$ and $r \in [0,1],$ in addition,
\begin{equation}
\lim_{t \to t_0}[F(t)]^r=\prod_{i=1}^n \Big[\lim_{t \to t_0} F_i^-(r,t),\lim_{t \to t_0} F_i^+(r,t) \Big],
\end{equation}
where $F_i^-(r,t)$ and $F_i^+(r,t)$ are real-valued functions.
\end{theorem}
\begin{proof}
Let $\lim_{t \to t_0} [F(t)]^r=[L]^r,$ for every $r \in [0,1].$ Then, $ \lim_{t \to t_0} d_L \big([F(t)]^r,[L]^r \big)=0.$ That is,
$$ \lim_{t \to t_0} \max_{1\leq i \leq n} \Big( \Big \{ \big |F^-_i(r,t)-L^-_i(r) \big|,\big |F^+_i(r,t)-L^+_i(r) \big | \Big \} \Big)=0.$$
As a result,
\begin{equation*}\label{34}
\max_{1\leq i \leq n} \Big( \Big \{ \Big |\lim_{t \to t_0} F^-_i(r,t)-L^-_i(r) \Big |,\Big |\lim_{t \to t_0} F^+_i(r,t)-L^+_i(r) \Big | \Big \} \Big)=0.
\end{equation*}
Hence, for any $r \in [0,1],$ we have
$$ \lim_{t \to t_0}F^-_i(r,t)=L^-_i(r) , \ \lim_{t \to t_0} F^+_i(r,t)=L^+_i(r).$$
Conversely, let
\begin{equation}\label{33}
\lim_{t \to t_0}F^-_i(r,t)=L^-_i(r) , \ \lim_{t \to t_0} F^+_i(r,t)=L^+_i(r).
\end{equation}
We show that $ \lim_{t \to t_0} [F(t)]^r=[L]^r \  \text{for} \ r \in [0,1].$
From \eqref{33}, we have
\begin{equation}\label{48}
 \max_{1\leq i \leq n} \Big \{ \big | F^-_i(r,t)-L^-_i(r) \big |,\big | F^+_i(r,t)-L^+_i(r) \big | \Big \} < \epsilon,
\end{equation}
for any $ i=1,2,...,n$ and $r \in [0,1].$ Consequently, $d_L ([F(t)]^r,[L]^r)<\epsilon. $
Which means
\begin{equation}\label{49}
 \lim_{t \to t_0} [F(t)]^r=[L]^r \  \text{for} \ r \in [0,1].
\end{equation}
Now, we show that $ \Big [\lim_{t \to t_0} F(t) \Big]^r=\lim_{t \to t_0} [F(t)]^r.$
By using Inequality \eqref{48},
\begin{equation}\label{50}
\max_{1\leq i \leq n} \Big( \Big \{ \Big |\lim_{t \to t_0} F^-_i(r,t)-L^-_i(r) \Big |,\Big |\lim_{t \to t_0} F^+_i(r,t)-L^+_i(r) \Big | \Big \} \Big)=0.
\end{equation}
Hence,
$$ d \Big( \Big[\lim_{t \to t_0} F(t) \Big]^r,[L]^r \Big)=0.$$
Therefore, according to Theorem \ref{z}(3), we have
$$ \Big [\lim_{t \to t_0} F(t) \Big]^r=[L]^r, \quad \text{for any} \  r \in [0,1].$$
Using the latter equality and \eqref{49}, we obtain
$$ \Big [\lim_{t \to t_0} F(t) \Big]^r=\lim_{t \to t_0} [F(t)]^r.$$
\end{proof}
In the following, we generalize the definitions of the directional derivative and the gradients proposed in \cite{Qiu} for a fuzzy $n$-cell number-valued function.
\begin{definition}\label{1}
Let $F: M \to L(E^n)$ be a fuzzy $n$-cell number-valued function and $t_0 \in M, d \in \mathbb{R}^m$ such that $t_0+td \in M$. If the following limit exists, then $ F'(t_0;d),$ the directional derivative of $F$ at $t_0$ in direction $d,$ defined as
$$ F'(t_0;d)=\lim_{t\rightarrow 0^+} \frac{F(t_0+td)\ominus_g F(t_0)}{t}.$$
\end{definition}
\begin{remark}\label{16}
In Definition \ref{1}, set $d=e_j,$ where $e_j$ is the $j$th canonical direction in $\mathbb{R}^m.$ If the following limit exists, then $F$ has partial derivative at the $j$th component:
$$ \frac{\partial F(t_0)}{\partial t_j^0}=\lim_{t\rightarrow 0} \frac{F(t_0+te_j)\ominus_g F(t_0)}{t},$$
 for any $j=1,2,...,m.$
\end{remark}
\begin{definition}\label{17}
Let $F: M \to L(E^n)$ be a fuzzy $n$-cell number-valued function. The gradient of $F$ at $t_0$ is a fuzzy vector as follows provided that all partial derivatives $F$ at $t_0$ exist:
$$\nabla F(t_0)=\Big(\frac{\partial F(t_0)}{\partial t_1^0},\frac{\partial F(t_0)}{\partial t_2^0},...,\frac{\partial F(t_0)}{\partial t_m^0}\Big).$$
\end{definition}
\begin{lemma}\label{3}
Let $F: M \to L(E^n)$ be a fuzzy $n$-cell number-valued function and $ t_0 \in M.$ If $\lim_{t \to t_0} [F(t)]^r=[u]^r$ for any $r \in [0,1]$ and  $u^-_i(r),u^+_i(r)$ satisfy Conditions in Theorem \ref{f}, then
$$ \lim_{t \to t_0} F(t)=u,$$ with $ [u]^r=\prod_{i=1}^n [u^-_i(r),u^+_i(r)].$
\end{lemma}
\begin{proof}
According to Definition \ref{18}, we must show that for any $\epsilon>0$ there exists a $ \delta>0$ such that for any $  t \in M,$ $||t-t_0|| < \delta $ implies
$$ D_L(F(t),u)< \epsilon.$$
Since $u^-_i(r),u^+_i(r)$ fulfill Conditions Theorem \ref{f}, thus, there exists a unique fuzzy n-cell number  $u \in L(E^n)$ such that $[u]^r=\prod_{i=1}^n [u^-_i(r),u^+_i(r)].$ Also, $\lim_{t \to t_0} [F(t)]^r=[u]^r$ implies that for any $\epsilon>0$ exists a $ \delta>0$ such that for any $  t \in M$ with $||t-t_0|| < \delta, $ we have $$d_L([F(t)]^r,[u]^r) < \epsilon, \quad \text{for any} \ r \in [0,1].$$
Therefore,
$$ \sup_{r\in [0,1]} d_L([F(t)]^r,[u]^r) \leq \epsilon,$$
that is, $$ D_L(F(t),u)\leq \epsilon.$$
\end{proof}
\begin{example}
Let $F:\mathbb{R} \to L(E) $ be a fuzzy $1$-cell number-valued function given by
\begin{equation*}
F(t)(x) = \left\{
\begin{matrix}
\frac{x}{|t|+1}+1, & -(|t|+1) \leq x\leq 0,\\
\frac{-x}{|t|+1}+1,  & 0 < x \leq (|t|+1),\\
0 & \text{otherwise}.
\end{matrix} \right.
\end{equation*}
We first determine the $r$-level sets of $F(t),$ \\
$$ [F(t)]^r=[(|t|+1)(r-1),(|t|+1)(1-r)].$$
In the following, we compute the directional derivative of $F$ at $t_0=0,$ in direction $d=1.$
\begin{align*}
\hspace{-1cm}\lim_{t\to 0^+} \Big[&\frac{F(t_0+td)\ominus_g F(t_0)}{t}\Big]^r \\
& \hspace*{1cm}=\lim_{t\to 0^+}\bigg[\frac{F(t)\ominus_g F(0)}{t}\bigg]^r\\
& \hspace*{1cm}=\lim_{t\to 0^+}\bigg[\inf_{\beta \geq r} \min \Big\{\frac{F^-(\beta,t)-F^-(\beta,0)}{t},\frac{F^+(\beta,t)-F^+(\beta,0)}{t} \Big \},\\
& \hspace*{2.5cm} \sup_{\beta \geq r} \max \Big \{\frac{F^-(\beta,t)-F^-(\beta,0)}{t},\frac{F^+(\beta,t)-F^+(\beta,0)}{t} \Big \}\bigg]\\
& \hspace*{1cm}=\bigg[\inf_{\beta \geq r} \min \Big\{ \lim_{t\to 0^+} \frac{F^-(\beta,t)-F^-(\beta,0)}{t},\lim_{t\to 0^+} \frac{F^+(\beta,t)-F^+(\beta,0)}{t} \Big \},\\
& \hspace{1.5cm}  \sup_{\beta \geq r} \max \Big \{ \lim_{t\to 0^+} \frac{F^-(\beta,t)-F^-(\beta,0)}{t},\lim_{t\to 0^+} \frac{F^+(\beta,t)-F^+(\beta,0)}{t} \Big \}\bigg]\\
& \hspace*{1cm}=\bigg[\inf_{\beta \geq r} \min \Big\{ \lim_{t\to 0^+} \frac{(|t|+1)(\beta-1)-(\beta-1)}{t},\lim_{t\to 0^+} \frac{(|t|+1)(1-\beta)-(1-\beta)}{t} \Big \},\\
& \hspace{1.5cm}  \sup_{\beta \geq r} \max \Big \{ \lim_{t\to 0^+} \frac{|t|(\beta-1)}{t},\lim_{t\to 0^+} \frac{|t|(1-\beta)}{t} \Big \}\bigg]\\
&\hspace*{1cm}=[r-1,1-r].
\end{align*}
Therefore,
$$ \lim_{t\rightarrow 0^+} \Big[\frac{F(t_0+td)\ominus_g F(t_0)}{t}\Big]^r=[r-1,1-r]=[u]^r,$$ for any $r\in[0,1].$
From Lemma \ref{3}, we have
$$ F'(0;1)=\lim_{t\rightarrow 0^+} \frac{F(0+t)\ominus_g F(0)}{t}=u(x).$$
Hence,
\begin{equation}
u(x) = \left\{
\begin{matrix}
x+1, & -1 \leq x \leq 0, \\
1-x, & 0 \leq x \leq 1,\\
0 & \text{otherwise}.
\end{matrix} \right.
\end{equation}
\end{example}
\begin{definition}\cite{Gong1} \label{46}
Let $F: M \to L(E^n)$ be a fuzzy $n$-cell number-valued function and $t_0, t \in \text{int} M.$ If $F(t) \ominus_g F(t_0)$ and a $u \in (L(E^n))^m$ exist such that
$$ \lim_{t \to t_0} \frac{D(F(t)\ominus_g F(t_0),u(t-t_0))}{||t-t_0||}=0,$$
then $F$ is differentiable at $t_0.$
\end{definition}
\begin{theorem}\label{v}
Let $u,v \in L(E^n).$ Then $u=v$ if and only if $[u]^r=[v]^r$ for any $r \in [0,1].$
\end{theorem}
\begin{proof}
If $u=v,$ then $u(x)=v(x)$ for any $x \in \mathbb{R}^m.$ Consequently,
$$ [u]^r=\{x\in \mathbb{R}^m: u(x)\geq r \}=\{x\in \mathbb{R}^m: v(x)\geq r\}=[v]^r,$$ for any $r \in (0,1].$ Also, If $r=0,$ then
$$ \overline{\{x\in \mathbb{R}^m: u(x)>0\}}=\overline{\{x\in \mathbb{R}^m: v(x)>0\}}.$$
So for any $r \in [0,1],$ we get $[u]^r=[v]^r.$ \\
Conversely, If $[u]^r=[v]^r,$ for any $ r\in [0,1],$ we show $u=v.$ \\
If $u\neq v,$ then there is an $x \in \mathbb{R}^m$ such that $u(x)\neq v(x).$ Without losing the generality, we assume that $ u(x) > v(x),$ then $ x \in [u]^{u(x)},$ and $ x \notin [v]^{u(x)},$ so $[u]^{u(x)}\neq [v]^{u(x)},$ which contradicts $[u]^r=[v]^r, \ r\in [0,1].$ Therefore, $u=v.$
\end{proof}
\begin{theorem}\cite{Gong}\label{47}
Let $F: M  \to L(E^n)$ be a fuzzy $n$-cell number-valued function, $ t \in M$ and $ d \in \mathbb{R}^m.$ Then $F$ is differentiable in the direction $d$ at $t$ if and only if for any $r \in [0,1],$ the functions $F^-_i(r,t), F^+_i(r,t):M \to \mathbb{R},$ $i=1,2,...,n,$ are differentiable in the direction $d$ at $t,$ and
$$ \inf_{ \beta \geq r} \min \{ {F_i^-}'(t,d,\beta), {F_i^+}'(t,d,\beta) \}, \quad   \sup_{ \beta \geq r} \max \{ {F_i^-}'(t,d,\beta), {F_i^+}'(t,d,\beta) \}$$
stisfy Conditions $(1)-(4)$ of Theorem \ref{f}, and 
$$ \Big [F' (t,d)\Big]^r=\prod_{i=1}^n \Big[\inf_{\beta \geq r} \min \Big \{ {F_i^-}'(t,d,\beta),{F_i^+}'(t,d,\beta) \Big \},
\sup_{\beta \geq r} \max \Big \{ {F_i^-}'(t,d,\beta),{F_i^+}'(t,d,\beta) \Big \} \Big].$$
\end{theorem}
In the following, we give the relation between directional derivative and gradient of a fuzzy $n$-cell number-valued function.
\begin{theorem}\label{t}
Let $F: M \to L(E^n)$ be a fuzzy $n$-cell number-valued function and $F$ be differentiable in the direction $d$ at $t \in M.$ Then for any $d \in \mathbb{R}^m$
$$F'(t;d)=\nabla F(t).d.$$
\end{theorem}
\begin{proof}
Since $F$ is differentiable in the direction $d$ at $t \in M,$ therefore according to Theorem \ref{12} and \ref{47} we have
\begin{align*}
 \Big [F'(t;d)\Big]^r \\
& \hspace*{-1.5cm}=\prod_{i=1}^n \Big[\inf_{\beta \geq r} \min \Big \{ {F_i^-}'(t,d,\beta),{F_i^+}'(t,d,\beta) \Big \},
\sup_{\beta \geq r} \max \Big \{ {F_i^-}'(t,d,\beta),{F_i^+}'(t,d,\beta) \Big \} \Big] \\
 & \hspace*{-1.5cm} =\prod_{i=1}^n \Big[\inf_{\beta \geq r} \min \{\nabla F_i^- (t,d,\beta),\nabla F_i^+ (t,d,\beta) \}, 
  \sup_{\beta \geq r} \max \{\nabla F_i^- (t,d,\beta),\nabla F_i^+ (t,d,\beta) \} \Big] \\
& \hspace*{-1.5cm}=\Big [\nabla F(t)d \Big]^r.
\end{align*}
It follows that based on Theorem \ref{v},
 $$ F'(t;d)=\nabla F(t)d. $$
\end{proof}
\begin{theorem}\label{tt}
Let $F: M \to L(E^n)$ be a convex fuzzy $n$-cell number-valued function, $M$ be a convex set and $t \in M.$ Then
$$ F(y) \ominus_g F(t) \geq F'(t;y-t),$$ for every $ y \in M.$
\end{theorem}
\begin{proof}
The convexity of $F$ implies
\begin{align}\label{kk}
F(t+\lambda(y-t))\ominus_g F(t) &= F((1-\lambda)t+\lambda y)\ominus_g F(t) \nonumber \\
& \leq ((1-\lambda)F(t)+\lambda F(y)) \ominus_g F(t)
\end{align}
for any $y,t \in M$ and $\lambda \in (0,1].$
According to Definition \ref{k}, for any $r \in [0,1],$
\begin{align*}
 [((1-\lambda)F(t)+\lambda F(y)) \ominus_g F(t)]^r \\
&\hspace*{-4cm}=\prod_{i=1}^n \Big[\inf_{\beta \geq r} \min \{(1-\lambda) F^-_i(\beta,t)+\lambda F^-_i(\beta,y)-F^-_i(\beta,t),\\
&\hspace*{-2cm} (1-\lambda)F^+_i(\beta,t)+\lambda F^+_i(\beta,y)-F^+_i(\beta,t) \},\\
& \hspace*{-3cm} \sup_{\beta\geq r} \max \{(1-\lambda)F^-_i(\beta,t)+\lambda F^-_i(\beta,y)-F^-_i(\beta,t),\\
& \hspace*{-2cm} (1-\lambda)F^+_i(\beta,t)+\lambda F^+_i(\beta,y)-F^+_i(\beta,t) \} \Big]\\
& \hspace*{-4cm}=\prod_{i=1}^n \Big[\inf_{\beta \geq r} \min \{\lambda (F^-_i(\beta,y)-F^-_i(\beta,t)),\lambda (F^+_i(\beta,y)-F^+_i(\beta,t))\},\\
& \hspace*{-2.5cm} \sup_{\beta\geq r} \max \{\lambda (F^-_i(\beta,y)-F^-_i(\beta,t)),\lambda (F^+_i(\beta,y)-F^+_i(\beta,t)) \} \Big]\\
& \hspace*{-4cm}=[\lambda (F(y) \ominus_g F(t))]^r.
\end{align*}
Applying Theorem \ref{v}, we have
\begin{equation}\label{14}
((1-\lambda)F(t)+\lambda F(y)) \ominus_g F(t)=\lambda (F(y) \ominus_g F(t)).
\end{equation}
Now, using Definition \ref{1}, Inequality \eqref{kk} and Equality \eqref{14}, we conclude that
\begin{align*}
F'(t;y-t)&= \lim_{\lambda \to 0^+} \frac{F(t+\lambda(y-t))\ominus_g F(t)}{\lambda}\\
&=\lim_{\lambda \to 0^+} \frac{F((1-\lambda)t+\lambda y)\ominus_g F(t)}{\lambda}\\
& \leq \lim_{\lambda \to 0^+} \frac{((1-\lambda)F(t)+\lambda F(y)) \ominus_g F(t)}{\lambda}\\
&=\lim_{\lambda \to 0^+} \frac{\lambda (F(y) \ominus_g F(t))}{\lambda}\\
&=F(y) \ominus_g F(t).
\end{align*}
Therefore,
$$ F(y) \ominus_g F(t) \geq F'(t;y-t),$$ for any $ y \in M.$
\end{proof}
\subsection{\bf Subdifferential of a Fuzzy $n$-cell Number-Valued Function}
In the following, we define the subdifferential of a fuzzy $n$-cell number-valued function considering the partial order introduced in Definition \ref{c}.
\begin{definition}\label{15}
Let $F: M  \to L(E^n)$ be a fuzzy $n$-cell number-valued function and $t \in M.$ If $F(z) \ominus_g F(t)$ and $ v \in (L(E^n))^m$ exist such that for any $ z \in M,$
$$ F(z) \ominus_g F(t)\geq v.(z-t),$$
then $ v$ is called the subgradient of $ F$ at $ t.$ The subdifferential of F at t is denoted by $\partial F(t)$ which is the set of subgradients of $F$ at $t.$
\end{definition}
\begin{example}\label{58}
Consider the fuzzy 1-cell number-valued function $ F:\mathbb{R} \to L(E)$ as
\begin{equation*}
F(t)(x) = \left\{
\begin{matrix}
\frac{x}{|t|+1}, & 0 \leq x \leq |t|+1, \\
2- \frac{x}{|t|+1}, & |t|+1 \leq x \leq 2(|t|+1), \\
0 & \text{otherwise}.
\end{matrix} \right.
\end{equation*}
It is easy to see that $ [F(t)]^r=[r (|t|+1),(2-r) (|t|+1)],$ for any $ r \in [0,1].$
First we show that $F$ is not differentiable at $t=0.$
\begin{align*}
\lim_{h \to 0} \Big [\frac{ F(0+h) \ominus_g F(0)}{h} \Big]^r \\
& \hspace*{-4cm}=\bigg [\inf_{\beta \geq r} \min \Big \{ \lim_{h \to 0} \frac{F^-(\beta,h)-F^-(\beta,0)}{h},\lim_{h \to 0} \frac{F^+(\beta,h)-F^+(\beta,0)}{h} \Big \}, \\
& \hspace*{-3cm} \sup_{\beta \geq r} \max \Big \{ \lim_{h \to 0} \frac{F^-(\beta,h)-F^-(\beta,0)}{h},\lim_{h \to 0} \frac{F^+(\beta,h)-F^+(\beta,0)}{h} \Big \} \bigg]\\
& \hspace*{-4cm}=\bigg [\inf_{\beta \geq r} \min \Big \{ \lim_{h \to 0} \frac{\beta |h|}{h},\lim_{h \to 0} \frac{(2-\beta) |h|}{h} \Big \},\sup_{\beta \geq r} \max \Big \{ \lim_{h \to 0} \frac{\beta |h|}{h},\lim_{h \to 0} \frac{(2-\beta) |h|}{h} \Big \} \bigg].
\end{align*}
To continue, we consider two cases:\\
Case I: If $ h \to 0^+,$ then
\begin{align*}
&  \bigg [\inf_{\beta \geq r} \min \Big \{ \lim_{h \to 0^+} \frac{\beta |h|}{h},\lim_{h \to 0^+} \frac{(2-\beta) |h|}{h}  \Big \},\sup_{\beta \geq r} \max \Big \{ \lim_{h \to 0^+} \frac{\beta |h|}{h},\lim_{h \to 0^+} \frac{(2-\beta) |h|}{h} \Big \} \bigg] \nonumber \\
& = \bigg [\inf_{\beta \geq r} \min \Big \{ \lim_{h \to 0^+} \frac{\beta h}{h},\lim_{h \to 0^+} \frac{(2-\beta) h}{h} \Big \},\sup_{\beta \geq r} \max \Big \{ \lim_{h \to 0^+} \frac{\beta h}{h},\lim_{h \to 0^+} \frac{(2-\beta) h}{h} \Big \} \bigg] \nonumber \\
& =[r,2-r].
\end{align*}
Case II: If $ h \to 0^-,$ then
\begin{align*}
&  \bigg [\inf_{\beta \geq r} \min \Big \{ \lim_{h \to 0^-} \frac{\beta |h|}{h},\lim_{h \to 0^-} \frac{(2-\beta) |h|}{h}  \},\sup_{\beta \geq r} \max \Big \{ \lim_{h \to 0^-} \frac{\beta |h|}{h},\lim_{h \to 0^-} \frac{(2-\beta) |h|}{h} \} \bigg] \nonumber \\
& = \bigg [\inf_{\beta \geq r} \min \Big \{ \lim_{h \to 0^-} \frac{-\beta h}{h},\lim_{h \to 0^-} \frac{-(2-\beta) h}{h}  \},\sup_{\beta \geq r} \max \Big \{ \lim_{h \to 0^-} \frac{-\beta h}{h},\lim_{h \to 0^-} \frac{-(2-\beta) h}{h} \} \bigg] \nonumber \\
& =[r-2,-r].
\end{align*}
Since the left and right derivatives are not equal, thus, $F$ is not differentiable at $t=0.$ \\
Now, we verify the subdifferentiability of $F $ at $t=0.$
By Definition \ref{15}, $ v \in L(E)$ is a subgradient of $F$ at $t=0,$ if for any $z \in \mathbb{R},$
$$ F(z) \ominus_g F(0)\geq v.(z-0),$$
which is equivalent to
$$ [F(z) \ominus_g F(0)]^r \geq [v.z]^r, \quad   \text{for any} \ r \in [0,1].$$
By using Theorem \ref{k}, we have
\begin{align*}
\Big [\inf_{\beta \geq r} \min \{ F^-(\beta,z)-F^-(\beta,0),F^+(\beta,z)-F^+(\beta,0) \},\\
& \hspace*{-7cm} \sup_{\beta \geq r} \max \{ F^-(\beta,z)-F^-(\beta,0),F^+(\beta,z)-F^+(\beta,0) \} \Big]\\
& \hspace*{-5cm} \geq [\min \{z v^-(r),z v^+(r) \},\max \{z v^-(r),z v^+(r) \}],
\end{align*}
for any $z \in \mathbb{R}.$ Therefore,
\begin{align*}
 \Big [\inf_{\beta \geq r} \min \{ \beta  |z|,(2-\beta)|z| \},\sup_{\beta \geq r} \max \{ \beta  |z|,(2-\beta)|z| \} \Big]\\
&\hspace*{-8cm}=\big[r|z|,(2-r)|z| \big] \geq [\min \{z v^-(r),z v^+(r) \},\max \{z v^-(r),z v^+(r) \}].
\end{align*}
We consider two cases:\\
Case I: If $z \geq 0,$ then
\begin{align*}
\big[rz,(2-r)z \big] \geq [z v^-(r),z v^+(r)].
\end{align*}
Now, by Definition \ref{c}, we have
\begin{equation*}
rz \geq z v^-(r), \quad (2-r)z \geq z v^+(r), \quad \text{for any} \ r \in [0,1].
\end{equation*}
Therefore,
\begin{equation}\label{41}
r \geq  v^-(r), \quad 2-r \geq  v^+(r), \quad \text{for any} \ r \in [0,1].
\end{equation}
Case II: If $z < 0,$ then
\begin{align*}
[-rz,(r-2)z] \geq [z v^+(r),z v^-(r)].
\end{align*}
Again,  by Definition \ref{c}, we have
\begin{equation*}
(r-2)z \geq z v^-(r), \quad -rz \geq z v^+(r), \quad \text{for any} \ r \in [0,1].
\end{equation*}
Therefore,
\begin{equation}\label{42}
r-2 \leq  v^-(r), \quad -r \leq  v^+(r), \quad \text{for any} \ r \in [0,1].
\end{equation}
Hence, by \eqref{41} and \eqref{42}, we get the subdifferential set of $F$ at $t=0,$   
$$ \partial F(0)= \{ v \in L(E): \ r-2 \leq v^-(r) \leq r, \ -r \leq v^+(r) \leq 2-r \}.$$
\begin{figure}[!htbp]
					\centering
					\includegraphics[scale=0.6]{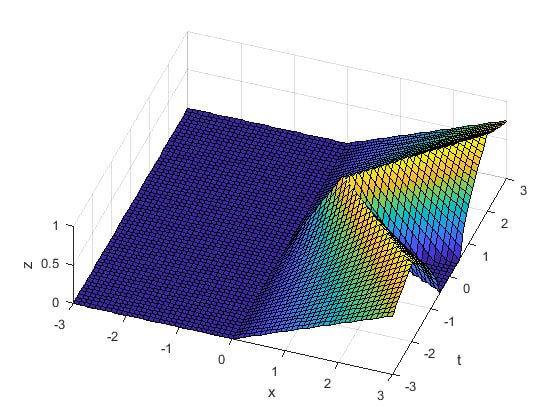}
					\caption{An example of a 1-cell fuzzy function which is not differentiable at zero.}
					\label{fig1}
				\end{figure}
\end{example}
\begin{definition}
A function $F: M \to L(E^n)$ is called subdifferentiable at $ t \in \text{int}(M)$ if $\partial F(t)\neq \emptyset.$
\end{definition}
\begin{theorem}
Let $F: M \to L(E^n)$ be a fuzzy $n$-cell number-valued function and $ t \in M.$ Then $\partial F(t)$ is a convex set in $L(E^n).$
\end{theorem}
\begin{proof}
If $\partial F(t)$ is empty or singleton, the result is obvious. Otherwise, for two subgradients $ u$ and $v$ that are arbitrary, we will show that
$$ \lambda u +(1-\lambda) v \in  \partial F(t), \ \text{for any} \  \lambda \in[0,1]. $$
Since $ u,v \in \partial F(t),$ then we have
$$ F(z) \ominus_g F(t) \geq (z-t).u,\qquad F(z) \ominus_g F(t) \geq (z-t).v, $$
for any $z \in M \subseteq \mathbb{R}^m,$ multiplying the first inequality by $\lambda,$ and the second inequality by $1-\lambda,$ and summing them and using Theorem \ref{xx}(5) yields
$$ F(z) \ominus_g F(t) \geq (z-t).(\lambda u+ (1-\lambda) v).$$
Therefore,
$$ \lambda u +(1-\lambda) v \in  \partial F(t),$$
That is, $ \partial F(t)$ is a convex set in $L(E^n)$.
\end{proof}
\begin{theorem}
Let $F: M \to L(E^n)$ be a fuzzy $n$-cell number-valued function and $ t \in M,$ then $\partial F(t)$ is a closed set.
\end{theorem}
\begin{proof}
Let $ \{v^{(k)} \}=\{(v_1^{(k)},v_2^{(k)},...,v_m^{(k)}) \} \in (L(E^n))^m$ be an arbitrary sequence of subgradient of $F$ at $t \in M$ such that
$$ \lim_{k \to \infty} \{v^{(k)} \}=v \in (L(E^n))^m.$$
We have
$$ \lim_{k \to \infty} \{v_j^{(k)} \}=v_j \in L(E^n), \ \text{for any} \ j=1,2,...,m.$$
We show that $ v=(v_1,v_2,...,v_m) \in \partial F(t).$ Since $ \lim_{k \to \infty} \{v_j^{(k)} \}=v_j,$ therefore, according to Theorem \ref{aa}, for any $\epsilon> 0$ there exists an $ N>0,$ such that
 $$ D_L(v_j^{(k)},v_j)<\epsilon,$$
for any $ k \geq N.$ i.e.,
$$\sup_{r \in [0,1]} d_L([v_j^{(k)}]^r,[v_j]^r)<\epsilon.$$
Since $ v_j^{(k)}, v_j \in L(E^n),$ So $ (v_{ij}^{(k)})^- (r), \ v_{ij}^- (r), \ (v_{ij}^{(k)})^+ (r), \ v_{ij}^+ (r) \in \mathbb{R},$ for any $j=1,2,...,m$ and $i=1,2,...,n.$ Therefore, by Theorem \ref{cc}, we have
$$ \sup_{r \in [0,1]} \max_{1 \leq i \leq n} \big \{ |(v_{ij}^{(k)})^- (r)-v_{ij}^- (r)|,|(v_{ij}^{(k)})^+ (r)-v_{ij}^+ (r)| \big \} <\epsilon, $$ for any $j=1,2,...,m.$ 
Thus,
$$ |(v_{ij}^{(k)})^-(r)-v_{ij}^-(r)|<\epsilon, \qquad |(v_{ij}^{(k)})^+(r)-v^+_{ij}(r)|<\epsilon$$\
 for any $r \in [0,1]$ and $j=1,2,...,m,$ which means that
 \begin{equation}\label{59}
 \lim_{k\to \infty} (v_{ij}^{(k)})^- (r)=v_{ij}^- (r), \qquad \lim_{k\to \infty} (v_{ij}^{(k)})^+ (r)=v_{ij}^+ (r)
 \end{equation}
 for any $r \in [0,1]$ and $ j=1,2,...,m.$\\
On the other hand, since $ \{v^{(k)} \} \in \partial F(t),$ we have
\begin{equation}\label{ff}
F(z) \ominus_g F(t)\geq v^{(k)} (z-t), \qquad \text{for any $z \in M.$}
\end{equation}
Using Definition \ref{c},
$$ \Big[F(z) \ominus_g F(t) \Big]^r \geq [v^{(k)} (z-t)]^r,$$
if and only if,
\begin{align*}
 \inf_{\beta \geq r} \min \{F^-_i(\beta,z)-F^-_i(\beta,t),F^+_i(\beta,z)-F^+_i(\beta,t)\}\\
& \hspace*{-4cm} \geq \min \bigg \{ \sum_{j=1}^m (v_{ij}^{(k)})^- (r) (z_j-t_j),\sum_{j=1}^m (v_{ij}^{(k)})^+(r) (z_j-t_j) \bigg \},
\end{align*}
and
\begin{align*}
\sup_{\beta \geq r} \max \{F^-_i(\beta,z)-F^-_i(\beta,t),F^+_i(\beta,z)-F^+_i(\beta,t)\}\\
 & \hspace*{-4cm} \geq \max \bigg \{ \sum_{j=1}^m (v_{ij}^{(k)})^- (r) (z_j-t_j),\sum_{j=1}^m (v_{ij}^{(k)})^+(r) (z_j-t_j) \bigg \}.
\end{align*}
In the last inequalities, taking $k \rightarrow \infty,$
\begin{align*}
 \inf_{\beta \geq r} \min \{F^-_i(\beta,z)-F^-_i(\beta,t),F^+_i(\beta,z)-F^+_i(\beta,t)\}\\
& \hspace*{-5cm} \geq \min \bigg \{ \lim_{k \rightarrow \infty} \sum_{j=1}^m (v_{ij}^{(k)})^- (r) (z_j-t_j), \lim_{k \rightarrow \infty} \sum_{j=1}^m (v_{ij}^{(k)})^+(r) (z_j-t_j) \bigg \},
\end{align*}
and
\begin{align*}
\sup_{\beta \geq r} \max \{F^-_i(\beta,z)-F^-_i(\beta,t),F^+_i(\beta,z)-F^+_i(\beta,t)\}\\
 & \hspace*{-5cm} \geq \max \bigg \{ \lim_{k \rightarrow \infty} \sum_{j=1}^m (v_{ij}^{(k)})^- (r) (z_j-t_j),\lim_{k \rightarrow \infty} \sum_{j=1}^m (v_{ij}^{(k)})^+(r) (z_j-t_j) \bigg \}.
\end{align*}
By Equalities \eqref{59}, we have 
\begin{align*}
 \inf_{\beta \geq r} \min \{F^-_i(\beta,z)-F^-_i(\beta,t),F^+_i(\beta,z)-F^+_i(\beta,t)\}\\
& \hspace*{-5cm} \geq \min \bigg \{  \sum_{j=1}^m v_{ij}^- (r) (z_j-t_j),  \sum_{j=1}^m v_{ij}^+(r) (z_j-t_j) \bigg \},
\end{align*}
and
\begin{align*}
\sup_{\beta \geq r} \max \{F^-_i(\beta,z)-F^-_i(\beta,t),F^+_i(\beta,z)-F^+_i(\beta,t)\}\\
 & \hspace*{-5cm} \geq \max \bigg \{  \sum_{j=1}^m v_{ij}^- (r) (z_j-t_j), \sum_{j=1}^m v_{ij}^+(r) (z_j-t_j) \bigg \}.
\end{align*}
Consequently, by Definition \ref{c}
$$ \Big[F(z) \ominus_g F(t) \Big]^r \geq [v (z-t)]^r $$
if and only if
$$ F(z) \ominus_g F(t)\geq v.(z-t), \qquad \text{for any $z \in M$}. $$
That is, $v \in (L(E^n))^m$ is a subdifferential of $F$ at $t \in M$ and hence $\partial F(t)$ is a closed set.
\end{proof}
\begin{lemma}\label{zz}
Let $u,v \in L(E^n)$ and $ u \ominus_g v $ exists, then
$$  u \ominus_g v \geq \hat{0} \qquad \text{if and only if} \qquad u \geq v.$$
\end{lemma}
\begin{proof}
Let $ u \ominus_g v \geq \hat{0}.$ Then applying Theorems \ref{bb}(1) and \ref{xx}(1), $u \ominus_g v  \geq v \ominus_g v,$ if and only if
$$(u \ominus_g v)+\hat{0}\geq (v \ominus_g v)+\hat{0} ,$$
 iff $$(u \ominus_g v)+(\hat{0}\ominus_g \hat{0}) \geq (v \ominus_g v)+ (\hat{0}\ominus_g \hat{0}) ,$$
iff (by Theorem \ref{mm}), $$(u + \hat{0}) \ominus_g (\hat{0}+v)\geq  (v + \hat{0}) \ominus_g(\hat{0} + v),$$
iff (using Theorem \ref{mm}), $$(u \ominus_g \hat{0})+(\hat{0} \ominus_g v)\geq (v \ominus_g \hat{0})+(\hat{0} \ominus_g v),$$
iff (using Theorem \ref{xx}(1)), $$u \ominus_g \hat{0} \geq v \ominus_g \hat{0} \ \text{if and only if} \ u \geq v.$$
\end{proof}
\begin{theorem}
Let $F: M \to L(E^n)$ be a fuzzy $n$-cell number-valued function and $M$ be a convex set. Assume that $\partial F(t) \neq \emptyset $ for $t \in M,$ then $F$ is convex.
\end{theorem}
\begin{proof}
Let $x,y \in M$ and $\lambda \in [0,1].$ Define $z=(1-\lambda)x+\lambda y.$ By the convexity of $M,$ $z \in M$ and $\partial F(z)\neq \emptyset. $ Hence there exists $ \xi \in \partial F(z)$ such that
$$F(x)\ominus_g F(z)\geq \xi.(x-z), \qquad F(y)\ominus_g F(z)\geq \xi.(y-z).$$ Therefore,
$$F(x)\ominus_g F(z)\geq \lambda \xi (x-y),\qquad F(y)\ominus_g F(z)\geq (1-\lambda) \xi.(y-x).$$
Multiplying the first inequality by $1-\lambda$ and the second by $\lambda,$
$$(1-\lambda)(F(x)\ominus_g F(z))\geq \lambda (1-\lambda) \xi (x-y),\qquad \lambda (F(y)\ominus_g F(z))\geq \lambda(1-\lambda) \xi (y-x).$$
Summing the last inequalities and using Theorem \ref{bb}(3), we obtain that
$$ \big((1-\lambda)F(x) \ominus_g (1-\lambda)F(z) \big)+(\lambda F(y) \ominus_g \lambda F(z))\geq \lambda (1-\lambda) \xi (x-y)+\lambda(1-\lambda) \xi (y-x).$$
By Theorem \ref{mm}, we get
$$ \big((1-\lambda)F(x)+\lambda F(y) \big) \ominus_g \big((1-\lambda)F(z)+\lambda F(z) \big)$$
$$ \qquad \geq \big(\lambda (1-\lambda)-\lambda (1-\lambda)\big)\xi (x-y),$$
and using Theorems \ref{bb}(4) and \ref{xx}(5), we obtain
$$  \big((1-\lambda)F(x)+\lambda F(y) \big) \ominus_g ((1- \lambda)+ \lambda) F(z) $$
$$ \hspace*{5cm} \geq \lambda (1-\lambda)\xi (x-y) \ominus_g \lambda (1-\lambda) \xi (x-y).$$
Also, applying Theorem \ref{bb}(1), we have
$$  \big((1-\lambda)F(x)+\lambda F(y) \big) \ominus_g F(z)\geq \hat{0}.$$
Therefore, using Lemma \ref{zz} we obtaion
$$(1-\lambda)F(x)+\lambda F(y) \geq F(z).$$
That is,
$$(1-\lambda)F(x)+\lambda F(y) \geq F((1-\lambda)x+\lambda y),$$
and $ F$ is convex.
\end{proof}
\begin{theorem}
Let $F: M \to L(E^n)$ be a convex fuzzy $n$-cell number-valued function. Then, the subdifferential of $F$ at $t \in M$ can be obtained by the subdifferential of the end points functions $F^-_i, F^+_i:M \to \mathbb{R},$ $i=1,2,...,n,$ and vice-versa.
\end{theorem}
\begin{proof}
Since $F$ is convex, therefore by Theorem \ref{52}, for any $r \in [0,1],$ $F^-_i(r,t)$ and $F^+_i(r,t),$ $i=1,2,...,n,$ are also convex. Hence, the functions $F^-_i(r,t)$ and $F^+_i(r,t)$ at $t \in M$ are subdifferentiable. Let $(v^-_{ij}(r))_{j=1}^m=(v^-_{i1}(r),v^-_{i2}(r),...,v^-_{im}(r)) \in \partial F^-_i(r,t)$ and $ (v^+_{ij}(r))_{j=1}^m=(v^+_{i1}(r),v^+_{i2}(r),...,v^+_{im}(r))  \in \partial F^+_i(r,t),$ $i=1,2,...,n.$ Then, we have
$$ F^-_i(r,z)-F^-_i(r,t) \geq  (v^-_{ij}(r))_{j=1}^m (z-t), \quad F^+_i(r,z)-F^+_i(r,t) \geq  (v^+_{ij}(r))_{j=1}^m (z-t),$$
for any $z \in M.$ Therefore,
$$ \inf_{\beta \geq r} \min \{F^-_i(\beta,z)-F^-_i(\beta,t),F^+_i(\beta,z)-F^+_i(\beta,t)\} \geq \min \bigg \{ \sum_{j=1}^m v^-_{ij} (r) (z_j-t_j),\sum_{j=1}^m v^+_{ij}(r) (z_j-t_j) \bigg \}$$
and,
$$ \sup_{\beta \geq r} \max \{F^-_i(\beta,z)-F^-_i(\beta,t),F^+_i(\beta,z)-F^+_i(\beta,t)\} \geq \max \bigg \{ \sum_{j=1}^m v^-_{ij} (r) (z_j-t_j),\sum_{j=1}^m v^+_{ij}(r) (z_j-t_j) \bigg \}.$$
Applying Definition \ref{c}, we have
$$ \Big[F(z) \ominus_g F(t) \Big]^r \geq [v (z-t)]^r, $$
and so,
$$ F(z) \ominus_g F(t)  \geq v (z-t), $$
for any $ z \in M.$ That is, $ v \in \partial F(t).$ 

Conversely, let $ v \in \partial F(t).$ Therefore,
\begin{equation}\label{53}
 F(z) \ominus_g F(t)  \geq v (z-t), 
\end{equation}
for any $ z \in M.$ We show that $F^-_i(r,t)$ and $F^+_i(r,t)$ at $t \in M$ are subdifferentiable. Using Definition \ref{c} and Inequality \eqref{53}, we have
$$ \inf_{\beta \geq r} \min \{F^-_i(\beta,z)-F^-_i(\beta,t),F^+_i(\beta,z)-F^+_i(\beta,t)\} \geq \min \bigg \{ \sum_{j=1}^m v^-_{ij} (r) (z_j-t_j),\sum_{j=1}^m v^+_{ij}(r) (z_j-t_j)  \bigg \}$$
and
$$ \sup_{\beta \geq r} \max \{F^-_i(\beta,z)-F^-_i(\beta,t),F^+_i(\beta,z)-F^+_i(\beta,t)\} \geq \max \bigg \{ \sum_{j=1}^m v^-_{ij} (r) (z_j-t_j),\sum_{j=1}^m v^+_{ij}(r) (z_j-t_j)  \bigg \}.$$
Now, Let $ \inf_{\beta \geq r} \min \{F^-_i(\beta,z)-F^-_i(\beta,t),F^+_i(\beta,z)-F^+_i(\beta,t)\}= F^-_i(r,z)-F^-_i(r,t)$
and $ \min \bigg \{ \sum_{j=1}^m v^-_{ij} (r) (z_j-t_j),\sum_{j=1}^m v^+_{ij}(r) (z_j-t_j)  \bigg \}=\sum_{j=1}^m v^-_{ij} (r) (z_j-t_j).$
In this case, 
\begin{equation}\label{54}
 F^-_i(r,z)-F^-_i(r,t) \geq \sum_{j=1}^m v^-_{ij} (r) (z_j-t_j), \ F^+_i(r,z)-F^+_i(r,t) \geq \sum_{j=1}^m v^+_{ij} (r) (z_j-t_j).
\end{equation}
That is, $  (v^-_{ij}(r))_{j=1}^m  \in \partial F^-_i(r,t)$ and $  (v^+_{ij}(r))_{j=1}^m  \in \partial F^+_i(r,t).$ Thus, $F^-_i(r,t), F^+_i(r,t):M \to \mathbb{R},$ $ i=1,2,...,n,$ are subdifferentiable at $ t \in M. $ 
\end{proof}
\begin{fact}
Let $F: M \to L(E^n)$ be a convex fuzzy $n$-cell number-valued function. Then, $F$ has a unique subgradient at $t \in M$ if and only if $F^-_i(r,t), F^+_i(r,t):M \to \mathbb{R},$ $i=1,2,...,n,$ have a unique subgradient at $t \in M.$ 
\end{fact}
\begin{theorem}\cite{Gong} \label{56}
Let $F: M \to L(E^n)$ be a fuzzy $n$-cell number-valued function and the functions $F^-_i(r,t), F^+_i(r,t):M \to \mathbb{R},$ $i=1,2,...,n,$ are differentiable at $t \in M,$ and
$$ \inf_{ \beta \geq r} \min \{ {F_{ij}^-}' (\beta,t),{F_{ij}^+}' (\beta,t) \}, \  \sup_{ \beta \geq r} \max \{ {F_{ij}^-}' (\beta,t),{F_{ij}^+}' (\beta,t) \} $$
satisfy Conditions $(1)-(4)$ of Theorem \ref{f}, where  ${F_{ij}^-}' (r,t)$ and ${F_{ij}^+}' (r,t)$ are the partial derivatives of $F^-_i(r,t)$ and $F^+_i(r,t)$ at $t \in M$ with respect to the $j$-th component, and   
$$ [F'_j(t)]^r=\prod_{i=1}^n \Big[ \inf_{ \beta \geq r} \min \{ {F_{ij}^-}' (\beta,t),{F_{ij}^+}' (\beta,t) \},\sup_{ \beta \geq r} \max \{ {F_{ij}^-}' (\beta,t),{F_{ij}^+}' (\beta,t) \} \Big].$$
Therefore, $F$ is differentiable at $t \in M.$
\end{theorem}
\begin{theorem}\label{57}
Let $F: M \to L(E^n)$ be a convex fuzzy $n$-cell number-valued function. If $F$ is differentiable at $t \in M,$ then
$$ \partial F(t)=\{\nabla F(t)\}.$$
Conversely, if $F$ has a unique subgradient at $t,$ then it is differentiable at $t$ and 
$$ \partial F(t)=\{\nabla F(t)\}.$$
\end{theorem}
\begin{proof}
Let $v \in \partial F(t),$ we have
$$ (y-t).v \leq  F(y) \ominus_g F(t),$$ for any $ y \in M.$
For any $ \lambda> 0 $ and $ d \in \mathbb{R}^m,$ set $ y=t + \lambda d.$ In this case,
$$ d.v \leq \frac{F(t + \lambda d) \ominus_g F(t)}{\lambda} .$$
Taking $ \lambda \rightarrow 0^+$ and using Theorem \ref{t}, we have
\begin{equation}
d.v \leq  \lim_{\lambda \to 0^+} \frac{F(t + \lambda d) \ominus_g F(t)}{\lambda}=F'(t;d)=\nabla F(t).d.
\end{equation}
Applying Theorem \ref{tt}, we obtain that
\begin{equation}\label{28}
d.v \leq \nabla F(t).d=\nabla F(t).(y-t)=F'(t;y-t) \leq F(y) \ominus_g F(t).
\end{equation}
That is, $ \nabla F(t) \in \partial F(t).$
Replacing $d$  by $-d$ in \eqref{28}, we get
\begin{equation*}
(-d).v \leq \nabla F(t).(-d),
\end{equation*}
and therefore,
\begin{equation}\label{29}
(-1).d.v \leq (-1).d.\nabla F(t).
\end{equation}
We get from \eqref{28} and \eqref{29},
$$ \sum_{i=1}^m d_i.v_i =d.v= \nabla F(t).d=\sum_{i=1}^m \frac{\partial F(t)}{\partial t_i} d_i . $$
Put $d=e_i $ for any $i=1,2,...,m.$ Therefore,
$$ v_i=\frac{\partial F(t)}{\partial t_i}, \  \text{for any} \  i=1,2,...,m. $$
Consequently,
$$ v= \nabla F(t).$$ i.e.,
$$ \partial F(t)=\{\nabla F(t)\}.$$

Conversely, let $ v \in (L(E^n))^m $ be a unique subgradient of $F$ at $ t \in M.$ Therefore, the left and right end-points of the $r$-level set of $v,$ i.e., 
$$ \inf_{\beta \geq r} \min \{v^-_{ij}(\beta),v^+_{ij}(\beta) \}, \  \sup_{\beta \geq r} \max \{v^-_{ij}(\beta),v^+_{ij}(\beta) \} $$ are unique subgradients of $F^-_i (r,t)$ and $F^+_i (r,t) $ at $ t \in M,$ respectively and satisfy Conditions (1)-(4) of Theorem \ref{f}, and we have
\begin{align*}
 \Big [ \partial F_j(t) \Big]^r &=\prod_{i=1}^n \Big[\inf_{\beta \geq r} \min \Big \{  \partial F^-_{ij} (\beta,t),\partial F^+_{ij} (\beta,t) \Big \}, \sup_{\beta \geq r} \max \Big \{ \partial F^-_{ij} (\beta,t),\partial F^+_{ij} (\beta,t) \Big \} \Big]\\
 & =\prod_{i=1}^n \Big[\inf_{\beta \geq r} \min \Big \{  v^-_{ij} (\beta), v^+_{ij} (\beta) \Big \}, \sup_{\beta \geq r} \max \Big \{ v^-_{ij} (\beta),v^+_{ij} (\beta) \Big \} \Big]\\
& =\prod_{i=1}^n \Big[\inf_{\beta \geq r} \min \Big \{  \nabla F^-_{ij} (\beta,t), \nabla F^+_{ij} (\beta,t) \Big \}, \sup_{\beta \geq r} \max \Big \{ \nabla F^-_{ij} (\beta,t),\nabla F^+_{ij} (\beta,t) \Big \} \Big]\\
&=\Big [ \nabla F_j(t) \Big]^r
\end{align*}
for any $ j=1,2,...,m.$ Based on Theorem \ref{v}, we have $ \partial F_j(t)=\nabla F_j(t),$ for any $ j=1,2,...,m.$
Note that, the uniqueness of the subgradients of $F_i^- (r,t)$ and $F_i^+ (r,t) $ implies the differentiability of  at $ t \in M.$ Since $F_i^- (r,t)$ and $F_i^+ (r,t) $ satisfy Conditions (1)-(4) of Theorem \ref{f}, Hence, $  {F_{ij}^-}' (r,t)= \lim_{h \to 0} \frac{F_i^- (r,t+h)-F_i^-(r,t)}{h}$ and $ {F_{ij}^+}' (r,t)= \lim_{h \to 0} \frac{F^+_i(r,t+h)-F^+_i(r,t)}{h}$ and as a result $ \inf_{ \beta \geq r} \min \{ {F_{ij}^-}' (\beta,t),{F_{ij}^+}' (\beta,t) \} $ and $ \sup_{ \beta \geq r} \max \{ {F_{ij}^-}' (\beta,t),{F_{ij}^+}' (\beta,t) \} $ Also satisfy Conditions (1)-(4) of Theorem \ref{f}. So, there exists a unique $ F'_j \in L(E^n)$ such that 
$$ [F'_j(t)]^r=\prod_{i=1}^n \Big[ \inf_{ \beta \geq r} \min \{ {F_{ij}^-}' (\beta,t),{F_{ij}^+}' (\beta,t) \},\sup_{ \beta \geq r} \max \{ {F_{ij}^-}' (\beta,t),{F_{ij}^+}' (\beta,t) \} \Big]$$ 
for any $ r \in [0,1], \ j=1,2,...,m.$ Now, we conclude that according to Theorem \ref{56}, $F$ is differentiable at $t \in M.$
\end{proof}
\begin{theorem}
Suppose $F,G:M \to L(E^n)$ are fuzzy $n$-cell number-valued functions and $\lambda >0.$ Then for any $t_0 \in M$
\begin{itemize}
\item[I)] $\partial(\lambda F)(t_0)=\lambda \partial F(t_0),$
\item[II)] $\partial (F+G)(t_0)=\partial F(t_0)+\partial G(t_0).$
\end{itemize}
\end{theorem}
\begin{proof}
\begin{itemize}
\item[I)] We show $ v \in \partial(\lambda F)(t_0)$  if and only if $ v \in \lambda \partial F(t_0).$ Using Theorems \ref{bb}(3) and \ref{xx}(2), for any $t \in M,$ we have
$$ v \in \partial(\lambda F)(t_0) \quad \text{iff} \quad (\lambda F)(t) \ominus_g (\lambda F)(t_0)\geq v.(t-t_0)$$
iff  $$\lambda (F(t) \ominus_g F(t_0))\geq v.(t-t_0)$$ 
 iff $$ F(t) \ominus_g F(t_0)\geq \frac{v}{\lambda} (t-t_0)$$
iff $$ \frac{v}{\lambda} \in \partial F(t_0).$$
Therefore,
$$\partial(\lambda F)(t_0)=\lambda \partial F(t_0) \qquad \text{for any} \quad \lambda >0.$$
\item[II)] We show that $ v \in \partial F(t_0)+\partial G(t_0)$ if and only if $ v \in \partial (F+G)(t_0).$ Let $ v \in \partial F(t_0)+\partial G(t_0).$ Then there exist $v_1 \in \partial F(t_0)$ and $v_2 \in \partial G(t_0)$ such that $v=v_1+v_2.$ By the definition, we have
$$F(t) \ominus_g F(t_0)\geq v_1.(t-t_0), \ G(t) \ominus_g G(t_0)\geq v_2.(t-t_0),$$
 for any $t \in M.$ Therefore,
$$(F(t) \ominus_g F(t_0))+(G(t) \ominus_g G(t_0))\geq v_1.(t-t_0)+ v_2.(t-t_0),$$
 for any $t \in M.$ Using Theorem \ref{mm}, we get
$$(F+G)(t) \ominus_g (F+G)(t_0)\geq (v_1+v_2).(t-t_0).$$
That is, $(v_1+v_2)=v \in \partial (F+G)(t_0).$ Therefore,
$$\partial F(t_0)+\partial G(t_0) \subseteq \partial (F+G)(t_0).$$
The converse is similarly obtained.
\end{itemize}
\end{proof}
Next, we give an example of Theorem \eqref{57}.
The matrix $A=(a_{kj})_{n \times n}, \ k,j=1,...,n,$ is a fuzzy matrix if its entries are fuzzy numbers, i.e., $ a_{kj}\in E.$
the operations on fuzzy matrices are defined by operations between fuzzy numbers.
\begin{definition}
An $ n \times n$ fuzzy matrix $A$ is said to be positive semidefinite and symmetric if $ A^-(r)=(a_{kj}^-(r))_{n \times n}$ and $ A^+(r)=(a_{kj}^+(r))_{n \times n}$ are $ n \times n$ positive semidefinite and symmetric for any $r \in [0,1],$ where $A^-(r)$ and $ A^+(r)$ are right and left end matrices $A.$
\end{definition}
\begin{lemma} \cite{Yuan}
If fuzzy matrix $A$ is positive semidefinite and symmetric, $ x \geq 0,$ then $x^T A x$ is a fuzzy number.
\end{lemma}
\begin{example}\label{62}
We show quadratic function $ F(x)=\frac{1}{2} x^T A x+b^T x =\frac{1}{2} \sum_{k=1}^n \sum_{j=1}^n a_{kj} x_k x_j + \sum_{j=1}^n b_j x_j $ has a unique subgradient at $x \in \mathbb{R}^n,$ where $A=(a_{kj})_{n \times n}$ is a positive semidefinite and symmetric matrix $(a_{kj}=a_{jk})$ and  $b=(b_1,b_2,...,b_n)^T.$ $ a_{kj}, \ b_j \in E.$ 
By using Definition \ref{17}, we first obtain the gradient of $F$ at $x_0.$ 
\begin{align*}
 \frac{\partial F(x_0)}{\partial x_j^0}\\
 & \hspace*{-1.5cm}=\lim_{t\rightarrow 0} \frac{F(x_0+te_j)\ominus_g F(x_0)}{t}\\
   &\hspace*{-1.5cm}= \lim_{t\rightarrow 0} \frac{\Big(\frac{1}{2} (x_0+te_j)^T A (x_0+te_j)+b^T (x_0+te_j)\Big)\ominus_g (\frac{1}{2} x_0^T A x_0+b^Tx_0)}{t}\\
  &\hspace*{-1.5cm}= \lim_{t\rightarrow 0} \frac{ \Big(\frac{1}{2}(x_0^T A x_0 +t(x_0^T A e_j)+t(e_j^T A x_0)+t^2 e_j^TAe_j)+b^T x_0 +tb^T e_j) \Big) \ominus_g (\frac{1}{2} x_0^T A x_0+b^Tx_0)}{t}\\
  &\hspace*{-1.5cm}= \lim_{t\rightarrow 0} \frac{\frac{1}{2}(t(x_0^T A e_j)+t(e_j^T A x_0)+t^2 e_j^TAe_j)+tb^T e_j}{t}=\frac{1}{2}\big((x_0^T A e_j)+(e_j^T A x_0) \big)+b^T e_j\\
&\hspace*{-1.5cm} =\frac{1}{2} \big((x_0^T A e_j)+(x_0^T A e_j)^T \big)+b^T e_j=x_0^T A e_j+b^T e_j= \sum_{k=1}^n  a_{kj} x^0_k+b_j \\
\end{align*}
Therefore,
$$ \nabla F(x_0)=\bigg(\frac{\partial F(x_0)}{\partial x_1^0},\frac{\partial F(x_0)}{\partial x_2^0},...,\frac{\partial F(x_0)}{\partial x_n^0}\bigg)=\bigg(\sum_{k=1}^n  a_{k1} x^0_k+b_1 ,...,\sum_{k=1}^n  a_{kn} x^0_k+b_n \bigg)=x_0^TA+b.$$
\end{example}
\subsection{\bf Optimality Conditions For Fuzzy Optimization Problems}
The topic of our discussion in this subsection is minimizing a fuzzy n-cell number-valued function and dealing with optimality conditions for unconstrained and constrained optimization problems. Consider a fuzzy function $F:M \to L(E^n).$ The following problem is called an unrestricted optimization problem
\begin{equation} \label{25}
\min_{x \in M} F(x).
\end{equation}
A feasible solution is defined as $x \in M.$
\begin{definition}
Let $F:M \to L(E^n)$ be a fuzzy $n$-cell number-valued function. A point $x_0 \in M$ is a global minimum of \eqref{25}, if $ F(x)\geq F(x_0)$ for any $x \in M.$ A point $x_0 \in M$ is a local minimum of \eqref{25}, if exists a $\delta>0$ such that $ F(x)\geq F(x_0)$ for any $x \in M\cap N(x_0,\delta)$ (where $ N(x_0,\delta)=\{ x \in \mathbb{R}^m: ||x-x_0|| < \delta \}).$
\end{definition}
\begin{theorem}\label{x}
Assume that $M$ is a convex set and $F:M \to L(E^n)$ is a convex fuzzy $n$-cell number-valued function. If $x_0 \in M$ is a local minimum of \eqref{25} then it is a global minimum of \eqref{25} in $M$ and conversely.
\end{theorem}
\begin{proof}
It is obvious that if $x_0$ is a global minimum solution of \eqref{25} in $M,$ then it is also a local minimum solution. Conversely, let $x_0$ be a local minimum solution. Since $x_0$ is a local minimum, there exists $\delta>0$ such that $ F(y)\geq F(x_0)$ for any $y \in M\cap N(x_0,\delta).$ Now, it suffices to show that for any $x \in M,$ $ F(x)\geq F(x_0).$ We consider two cases.\\\\
Case I. Let $||x-x_0||<\delta.$ In thise case, $ x \in N(x_0,\delta)$ and therefore, $ F(x)\geq F(x_0).$\\\\
Case II. Let $||x-x_0|| \geq\delta.$ Then, $0 \leq 1-\frac{\delta}{||x-x_0||} < 1.$ Take $ 1-\frac{\delta}{||x-x_0||}<k<1,$ then $(1-k){||x-x_0||}<\delta.$ Therefore,
$$ ||(kx_0+(1-k)x)-x_0||=(1-k)||x-x_0||<\delta.$$
We conclude $kx_0+(1-k)x \in N(x_0,\delta).$ The convexity of $M$ implies $kx_0+(1-k)x \in N(x_0,\delta)\cap M,$ and therefore, $F(kx_0+(1-k)x)\geq F(x_0).$ Also, by using the convexity of $F,$ we have
$$ F(x_0) \leq k F(x_0)+ (1-k) F(x),$$
if and only if by Definition \ref{c}, $i=1,2,...,n,$ for any $r \in [0,1],$
\begin{equation} \label{50}
 F^-_i(r,x_0) \leq k F_i^-(r,x_0)+ (1-k) F_i^-(r,x),
\end{equation}
\begin{equation} \label{51}
F^+_i(r,x_0) \leq k F_i^+(r,x_0)+ (1-k) F_i^+(r,x).
\end{equation}
Therefore, from \eqref{50} and \eqref{51},
\begin{align*}
(1-k) F^-_i(r,x_0) \leq (1-k) F_i^-(r,x),\\
(1-k) F^+_i(r,x_0) \leq  (1-k) F_i^+(r,x).
\end{align*}
By Definition \ref{c},
$$ (1-k) F(x_0) \leq (1-k) F(x).$$
That is,
$$ F(x_0)\leq F(x),$$ and $x_0 \in M$ is a global minimum solution.
\end{proof}
\begin{theorem}\label{24}
Let $M$ be a convex set and $F:M \to L(E^n)$ be a convex fuzzy $n$-cell number-valued function. Then $ x_0=(x^0_1,x^0_2,...,x^0_m) $ is a minimum solution of \eqref{25} if and only if $(\hat{0},\hat{0},...,\hat{0}) \in \partial F(x_0).$
\end{theorem}
\begin{proof}
By definition, $u={\bf \hat{0}}=(\hat{0},\hat{0},...,\hat{0}) \in \partial F(x_0)$ if and only if for any $x \in M,$
$$F(x)\ominus_g F(x_0) \geq (x-x^0).u =(x-x^0).{\bf {\hat{0}}}={\bf {\hat{0}}},$$
if and only if by Lemma \ref{zz},
$$F(x)\geq F(x_0).$$
That is, $ x_0 $ is a minimum solution of \eqref{25}.
\end{proof}
\begin{example}
In Example \ref{62}, suppose $A= \left(
\begin{array}{ccc}
\hat{a} & \hat{0}  \\
\hat{0} & \hat{a}
\end{array} \right),$ $b= \left(
\begin{array}{ccc}
\hat{0} & \hat{0} \\
\end{array} \right),$
where $\hat{a}=(0,4,8)$ is a triangular fuzzy number and $x=(x_1,x_2) \in \mathbb{R}^2_{+}.$ The quadratic function can be obtained as follows:
\begin{align*}
F(x) &=\frac{1}{2} x^T A x+b^T x =\frac{1}{2} (x_1 \ x_2)  \left(
\begin{array}{ccc}
\hat{a} & \hat{0}  \\
\hat{0} & \hat{a}
\end{array} \right) \left(
\begin{array}{ccc}
x_1  \\
x_2
\end{array} \right)+  \left(
\begin{array}{ccc}
\hat{0} & \hat{0} \\
\end{array} \right)\left(
\begin{array}{ccc}
x_1  \\
x_2
\end{array} \right)\\
 &=(0,2,4) x_1^2 + (0,2,4) x_2^2.
\end{align*}
The function of $F$ obtains the minimum at $x=(0,0).$ According to Example \ref{62}, $F$ is differentiable at $x,$ and therefore By Theorem \ref{57}, we have $ \partial F(x)=\{\nabla F(x)\}.$ Also, 
$$ \nabla F(x)= Ax+b^T=\left(
\begin{array}{ccc}
\hat{a} & \hat{0}  \\
\hat{0} & \hat{a}
\end{array} \right) \left(
\begin{array}{ccc}
x_1   \\
x_2
\end{array} \right) 
+  \left(
\begin{array}{ccc}
\hat{0} \\
\hat{0}
\end{array} \right)= \hat{a} x_1 + \hat{a} x_2+\hat{0}.$$
Hence, $\nabla F(0,0)=\hat{0}.$ That is, $ \hat{0} \in  \partial F(0,0).$
\end{example}
\subsubsection{\bf The Lagrangian dual problem}
Consider the following fuzzy constrained optimization problem,
\begin{equation}\label{pp}
\min\{F(x): G(x)\leq \hat{0}, \ x \in M\},
\end{equation}
where $F:M \to L(E^n)$ is a convex fuzzy $n$-cell number-valued function and
$G:M \to (L(E^n))^k$ is a convex $k$-dimensional fuzzy $n$-cell vector-valued function.
The Lagrangian $L:\mathbb{R}^m \times \mathbb{R}^k \to L(E^n)$ associated with \eqref{pp} is
$$L(x,\lambda)=F(x)+ \lambda.G(x)=F(x)+\sum_{j=1}^K \lambda_j G_j(x),$$
where ${\bf \lambda}=(\lambda_1,...,\lambda_k)$ and $ G(x)=(G_1(x),...,G_k(x)),$ $\lambda_j\geq0 \ (j=1,2,...,k).$
The Lagrangian dual problem is presented below
$$ \max_{\lambda\geq 0}d(\lambda),$$
where $ d(\lambda)=\min_{x \in M} L(x,\lambda).$
\begin{theorem}(Karush–Kuhn–Tucker optimality conditions)\label{26}
Let $x^*$ be an optimal solution of \eqref{pp} and exists $\hat{x} \in M$ with $G_j (\hat{x}) < \hat{0},$ $j=1,2,...,k.$ Then there exist $ \lambda_1, \lambda_2,...,\lambda_k \geq 0$ such that
\begin{equation}\label{20}
\hat{0} \in \partial F(x^*) + \lambda. \partial G(x^*)=\partial F(x^*) + \sum_{j=1}^k  \lambda_j. \partial G_j(x^*)
\end{equation}
\begin{equation}\label{21}
\lambda_j G_j(x^*) =\hat{0}, \ j=1,2,...,k.
\end{equation}
Conversely, $ x^*$ is an optimal solution of \eqref{pp}, if it fulfills requirements of \eqref{20} and \eqref{21} for certain $ \lambda_1,\lambda_2,...,\lambda_k  \geq 0.$
\end{theorem}
\begin{proof}
Suppose $x^*$ is an optimal solution of \eqref{pp}. Therefore, for any $x \in M,$
$$ F(x^*) \leq F(x), \ G_j(x^*) \leq \hat{0}, \ j=1,2,...,k. $$
Based on Definition \ref{c}, we have
$$ F^-_i(r,x^*) \leq F^-_i(r,x), \quad G^-_{ij}(r,x^*) \leq 0 $$ and
$$ F^+_i(r,x^*) \leq F^+_i(r,x), \quad G^+_{ij}(r,x^*) \leq 0,$$
for any $i=1,2,...,n $ and $r \in [0,1].$ Since $ F^+_i(r,x), \ F^-_i(r,x),$ $G^-_{ij}(r,x^*)$ and $G^+_{ij}(r,x^*)$ are convex real-valued functions with respect to $r,$ thus for every $r \in [0,1],$
\begin{equation}\label{22}
 0 \in  \partial F^-_i(r,x^*)+ \sum_{j=1}^k \lambda_j \partial G^-_{ij}(r,x^*), \quad \lambda_j G^-_{ij}(r,x^*)=0
\end{equation}
and
\begin{equation}\label{23}
 0 \in \partial F^+_i(r,x^*)+ \sum_{j=1}^k \lambda_j \partial G^+_{ij}(r,x^*), \quad \lambda_j G^+_{ij}(r,x^*)=0.
\end{equation}
Therefore, from the definition of subgradient of real-valued functions, we have
$$ \Big (F^-_i(r,z) + \sum_{j=1}^k \lambda_j  G^-_{ij}(r,z) \Big)- \Big (F^-_i(r,x^*) + \sum_{j=1}^k \lambda_j  G^-_{ij}(r,x^*) \Big) \geq 0.(z-x^*)=0$$ and
$$ \Big(F^+_i(r,z) + \sum_{j=1}^k \lambda_j  G^+_{ij}(r,z) \Big)- \Big(F^+_i(r,x^*) + \sum_{j=1}^k \lambda_j  G^+_{ij}(r,x^*) \Big) \geq 0.(z-x^*)=0,$$
for any $z \in M.$ Therefore,
\begin{align*}
\inf_{\beta \geq r} \min \Big \{ \Big ( F^-_i(\beta,z) + \sum_{j=1}^k \lambda_j  G^-_{ij}(\beta,z) \Big)- \Big (F^-_i(\beta,x^*) + \sum_{j=1}^k \lambda_j  G^-_{ij}(\beta,x^*) \Big),\\
& \hspace*{-10cm} \Big(F^+_i(\beta,z) + \sum_{j=1}^k \lambda_j  G^+_{ij}(\beta,z) \Big)- \Big (F^+_i(\beta,x^*) + \sum_{j=1}^k \lambda_j  G^+_{ij}(\beta,x^*) \Big) \Big \} \geq 0
\end{align*}
and
\begin{align*}
\sup_{\beta \geq r} \max \Big \{ \Big ( F^-_i(\beta,z) + \sum_{j=1}^k \lambda_j  G^-_{ij}(\beta,z) \Big)- \Big (F^-_i(\beta,x^*) + \sum_{j=1}^k \lambda_j  G^-_{ij}(\beta,x^*) \Big),\\
& \hspace*{-10cm} \Big(F^+_i(\beta,z) + \sum_{j=1}^k \lambda_j  G^+_{ij}(\beta,z) \Big)- \Big (F^+_i(\beta,x^*) + \sum_{j=1}^k \lambda_j  G^+_{ij}(\beta,x^*) \Big) \Big \} \geq 0.
\end{align*}
By Definition \ref{c}, we have
\begin{align*}
\prod_{i=1}^n \bigg[\inf_{\beta \geq r} \min \bigg \{ \Big ( F^-_i(\beta,z) + \sum_{j=1}^k \lambda_j  G^-_{ij}(\beta,z) \Big)- \Big (F^-_i(\beta,x^*) + \sum_{j=1}^k \lambda_j  G^-_{ij}(\beta,x^*) \Big),\\
& \hspace*{-10.5cm} \Big(F^+_i(\beta,z) + \sum_{j=1}^k \lambda_j  G^+_{ij}(\beta,z) \Big)- \Big (F^+_i(\beta,x^*) + \sum_{j=1}^k \lambda_j  G^+_{ij}(\beta,x^*) \Big) \bigg \},\\
\sup_{\beta \geq r} \max \bigg \{ \Big ( F^-_i(\beta,z) + \sum_{j=1}^k \lambda_j  G^-_{ij}(\beta,z) \Big)- \Big (F^-_i(\beta,x^*) + \sum_{j=1}^k \lambda_j  G^-_{ij}(\beta,x^*) \Big),\\
& \hspace*{-10.5cm} \Big(F^+_i(\beta,z) + \sum_{j=1}^k \lambda_j  G^+_{ij}(\beta,z) \Big)- \Big (F^+_i(\beta,x^*) + \sum_{j=1}^k \lambda_j  G^+_{ij}(\beta,x^*) \Big) \bigg \} \bigg]\\
 \geq 0.(z-x^*)=[\hat{0}]^r.
\end{align*}
Now, using Theorem \ref{k},
$$ \Big [\Big(F(z)+ \sum_{j=1}^k \lambda_j G_j(z) \Big) \ominus_g \Big(F(x^*)+ \sum_{j=1}^k \lambda_j G_j(x^*) \Big) \Big]^r \geq [\hat{0}]^r,$$
and therefore,
$$ \Big(F(z)+ \sum_{j=1}^k \lambda_j G_j(z) \Big) \ominus_g \Big(F(x^*)+ \sum_{j=1}^k \lambda_j G_j(x^*) \Big) \geq \hat{0},$$
for any $z \in M.$
It can be concluded that
$$ \hat{0} \in \partial F(x^*) + \lambda. \partial G(x^*).$$
Also, using Definition \ref{c} and Equalities \eqref{22} and \eqref{23}, we obtain that
$$ \lambda_j G_j(x^*) =\hat{0}, \ j=1,2,...,k.$$

Conversely, let $ x^*$ satisfies Conditions \eqref{20} and \eqref{21} for some $ \lambda_1,\lambda_2,...,\lambda_k  \geq 0.$ Suppose $x_0$ is a feasible solution of \eqref{pp}, that is, $ G_j(x_0) \leq \hat{0},$ for any $ j=1,2,...,k.$ It must be shown that $ F(x^*) \leq F(x_0).$ Put
$$ H(x):= F(x) + \sum_{j=1}^k \lambda_j G_j(x).$$
The convexity of $H(x)$ and \eqref{20} implies that
$$\hat{0} \in \partial H(x^*). $$
So by Theorem \ref{24}, $x^*$ is a minimum solution of $H$ on $M.$ Utilizing the Equality \eqref{21}, we have
$$ F(x^*) =F(x^*) + \sum_{j=1}^k \lambda_j G_j(x^*)=H(x^*) \leq H(x_0)= F(x_0) + \sum_{j=1}^k \lambda_j G_j(x_0)  \leq F(x_0).$$
This implies that $ F(x^*) \leq F(x_0),$ i.e., $x^*$ is an optimal solution of \eqref{pp}.
\end{proof}
\begin{example}
For given $t \in \mathbb{R},$ define the fuzzy functions $F: \mathbb{R}\to L(E)$ and $G:\mathbb{R} \to L(E^2)$ as
\begin{equation*}
F(t)(x)= \left\{
\begin{matrix}
1+x-t^2, &   x \in [t^2-1,t^2]\\
1-x+t^2, &  x \in [t^2,t^2+1]\\
\hspace*{-1.5cm} 0, & \hspace*{0.6cm} x \notin [t^2-1,t^2+1], \\
\end{matrix} \right.
\end{equation*}
and
\begin{equation*}
\hspace*{1cm} G(t)(x,y)= \left\{
\begin{matrix}
1+t+x, &   x \in [-t-1,-t], y \in [-2,-1]\\
1-t-x, &  x \in [-t,-t+1], y \in [-2,-1]\\
 0, & \text{otherwise} \\
\end{matrix} \right.
\end{equation*}
Then
$$[F(t)]^r=[t^2-1+r,t^2+1-r], \ [G(t)]^r=[-t-1+r,-t+1-r] \times [-2,-1]$$ for any $ r \in [0,1].$ For $t=2$ we have $ G^+_1(r,2),G^-_1(r,2) < 0.$ Therefore, exists $\hat{x} \in M$ with $G_j (\hat{x}) < \hat{0},$ $j=1,2.$
Since $ F(0) \leq F(t),$ for any $ t \in \mathbb{R},$ therefore $t_0=0$ is an optimal solution of \eqref{pp}. By the assumption of Theorem \ref{26}, we show that there exist $ \lambda_1, \lambda_2 \geq 0$ such that
\begin{equation}\label{60}
\hat{0} \in \partial F(0) + \lambda. \partial G(0)=\partial F(0) + \sum_{j=1}^2  \lambda_j. \partial G_j(0)
\end{equation}
and
\begin{equation}\label{61}
\lambda_j G_j(0) =\hat{0}, \ j=1,2.
\end{equation}
Let $ v \in \partial \big(F(0)+\sum_{j=1}^2  \lambda_j.  G_j(0) \big).$ By the definition \ref{15}, we show $ v=\hat{0},$ 
$$ \big(F(t) + \lambda_1 G_1(t) + \lambda_2 G_2(t) \big) \ominus_g \big (F(0) + \lambda_1 G_1(0) + \lambda_2 G_2(0) \big) \geq v.(t-0). $$
Applying Definition \ref{c}, we get
\begin{align*}
& \min \{t v^-(r),t v^+(r) \} \\
& \hspace*{-0.5cm} \leq \inf_{\beta \geq r} \min \big \{ F^-(\beta,t)+\lambda_1 G^-_1(\beta,t) + \lambda_2 G^-_2(\beta,t) -F^-(\beta,0)- \lambda_1 G^-_1(\beta,0) - \lambda_2 G^-_2(\beta,0) , \\
& \hspace*{1.5cm} F^+(\beta,t)+ \lambda_1 G^+_1(\beta,t) + \lambda_2 G^+_2(\beta,t) -F^+ (\beta,0)- \lambda_1 G^+_1(\beta,0) - \lambda_2 G^+_2(\beta,0) \big \} 
\end{align*}
and
\begin{align*}
& \max \{t v^-(r),t v^+(r) \} \\
& \hspace*{-0.5cm} \leq \sup_{\beta \geq r} \max \big \{ F^-(\beta,t)+ \lambda_1 G^-_1(\beta,t) + \lambda_2 G^-_2(\beta,t) -F^-(\beta,0)- \lambda_1 G^-_1(\beta,0) - \lambda_2 G^-_2(\beta,0) , \\
& \hspace*{1.5cm} F^+(\beta,t)+ \lambda_1 G^+_1(\beta,t) + \lambda_2 G^+_2(\beta,t) -F^+ (\beta,0)- \lambda_1 G^+_1(\beta,0) - \lambda_2 G^+_2(\beta,0) \big \}.
\end{align*}
Therefore,
\begin{align*}
& \min \{t v^-(r),t v^+(r) \} \\
& \leq \inf_{\beta \geq r} \min \big \{ t^2-1+\beta + \lambda_1 (-t-1+\beta)- 2 \lambda_2+1-\beta -\lambda_1(-1+\beta)+2 \lambda_2,\\
& \hspace*{1cm} t^2 +1-\beta+ \lambda_1 (-t +1-\beta)- \lambda_2 -1+\beta -\lambda_1 (1-\beta) +\lambda_2 \}
\end{align*}
and
\begin{align*}
& \max \{t v^-(r),t v^+(r) \} \\
& \leq \sup_{\beta \geq r} \max \big \{ t^2-1+\beta + \lambda_1 (-t-1+\beta)- 2 \lambda_2+1-\beta -\lambda_1(-1+\beta)+2 \lambda_2,\\
& \hspace*{1cm} t^2 +1-\beta+ \lambda_1 (-t +1-\beta)- \lambda_2 -1+\beta -\lambda_1 (1-\beta) +\lambda_2 \}.
\end{align*}
Hence,
\begin{align*}
& \hspace*{-1.2cm} \Big [ \min \{t v^-(r),t v^+(r) \},\max \{t v^-(r),t v^+(r) \} \Big ] \\
 & \leq \Big [ \inf_{\beta \geq r} \min \big \{ t^2 - \lambda_1 t ,t^2 - \lambda_1 t  \},   \sup_{\beta \geq r} \max \big \{ t^2 - \lambda_1 t ,t^2 - \lambda_1 t   \} \Big].  
\end{align*}
Thus,
$$ \Big [ \min \{t v^-(r),t v^+(r) \},\max \{t v^-(r),t v^+(r) \} \Big ] \leq \Big [ t^2 - \lambda_1 t ,t^2 - \lambda_1 t  \Big].$$
We consider two cases:\\
Case I: If $t \geq 0,$ therefore
$$ [t v^-(r), tv^+(r)] \leq \Big [ t^2 - \lambda_1 t ,t^2 - \lambda_1 t  \Big], $$
and so,
\begin{equation}\label{43}
[v^-(r), v^+(r)] \leq \Big [ t - \lambda_1  ,t - \lambda_1  \Big].
\end{equation}
Case II: If $t < 0,$ therefore
$$ [t v^+(r), tv^-(r)] \leq \Big [ t^2 - \lambda_1 t ,t^2 - \lambda_1 t  \Big], $$
and so,
\begin{equation}\label{44}
 [v^+(r), v^-(r)] \geq \Big [ t - \lambda_1  ,t - \lambda_1  \Big].
\end{equation}
From \eqref{43} and \eqref{44}, we have
$$  v^-(r)= t - \lambda_1, \quad v^+(r)= t - \lambda_1. $$
Therefore,
\begin{equation*}
v(x)=\widehat{(t-\lambda_1)}(x)= \left\{
\begin{matrix}
1 &   x =t-\lambda_1 \\
0 & x \neq t-\lambda_1.
\end{matrix} \right.
\end{equation*}
For any $ \lambda_1=t \geq 0,$ we have
\begin{equation*}
v(x)=\hat{0}(x)= \left\{
\begin{matrix}
1 &   x =0 \\
0 & x \neq 0
\end{matrix} \right.
\end{equation*}
That is,
$$ \hat{0} \in \partial F(0) + \sum_{j=1}^2  \lambda_j. \partial G_j(0).$$
Also, for $ \lambda_1=t \geq 0$ and $ \lambda_2 =0,$ the following equalities hold,
$$ \lambda_1 G^-_1(r,0)=\lambda_1 (-1+r)=0, \quad \lambda_2 G^-_2(r,0)=-2 \lambda_2 =0 $$ and
$$ \lambda_1  G^+_1(r,0)=\lambda_1 (1-r)=0, \quad  \lambda_2  G^+_2(r,0)=-\lambda_2 =0.$$
 Therefore, $ \lambda_j G_j(0) =\hat{0}, \ j=1,2.$ So, we show \eqref{60} and \eqref{61}.
\end{example}
\subsection{Optimality conditions for the composite problem}
In this subsection, we state the necessary and sufficient conditions for the optimal solution of a composite function by subgradient.
\begin{fact}\label{30}
Let $ u , v,w,z \in L(E^n).$ Then
$$ u+v \leq w+z  \quad \text{if and only if} \quad u \ominus_g w \leq z \ominus_g v.$$
\end{fact}
\begin{theorem}
Let $F,G:M \to L(E^n)$ such that $ \mathrm{dom(G)} \subseteq \mathrm{int(dom(F))}$ and $G$ be a convex function. Consider the problem
$$ (P) \quad \min_{x \in M} F(x)+G(x).$$
\begin{itemize}
\item[(I)]
(necessary condition) If $ x^* \in \mathrm{dom(G)}$ is a local optimal solution of $(P)$ and $F$ is differentiable at $x^*,$ then
\begin{equation}
-\nabla F(x^*) \in \partial G(x^*). \label{AN}
\end{equation}
\item[(II)]
(necessary and sufficient condition for convex problems). Assume that $F$ is convex fuzzy mapping. If $F$ is differentiable at $ x^* \in \mathrm{dom(G)}$, then $ x^*$ is a global optimal solution of $(P)$ if and only if \eqref{AN} is satisfied.
\end{itemize}
\end{theorem}
\begin{proof}
(I) Let $y \in int(M).$ Since $M$ is a convex set, then for any $ \lambda \in (0,1),$ we have $ x=(1-\lambda) x^*+\lambda y \in M.$ By using the local optimality of $ x^*,$ there exists $ \delta>0$ such that
$$ F(x^*)+G(x^*)\leq F(x)+G(x),$$ for any $ x \in M\cap N(x^*,\delta).$
That is,
$$ F(x^*)+G(x^*) \leq F((1-\lambda) x^*+\lambda y)+G((1-\lambda) x^*+\lambda y).$$
Since $G$ is a convex function,
$$F(x^*)+G(x^*) \leq F((1-\lambda) x^*+\lambda y)+(1-\lambda)G(x^*)+\lambda G(y).$$
Using Fact \eqref{30},
\begin{align*}
 F(x^*) \ominus_g F((1-\lambda) x^*+\lambda y) \leq \big((1-\lambda)G(x^*)+\lambda G(y) \big)\ominus_g G(x^*) \\
& \hspace*{-6.5cm}=\big(G(x^*)+(-\lambda) G(x^*) + \lambda G(y) \big) \ominus_g G(x^*)\\
&\hspace*{-6.5cm}= \big((G(x^*) \ominus_g \lambda G(x^*))+\lambda G(y)\big) \ominus_g G(x^*)\\
&\hspace*{-6.5cm}= \big ((G(x^*) \ominus_g \lambda G(x^*))+(\lambda G(y) \ominus_g \hat{0})\big) \ominus_g G(x^*)\\
&\hspace*{-6.5cm}= \big ((G(x^*) + \lambda G(y)) \ominus_g (\lambda G(x^*) + \hat{0}) \big) \ominus_g G(x^*)\\
&\hspace*{-6.5cm}= \big ((G(x^*) \ominus_g \hat{0} )+ (\lambda G(y) \ominus_g \lambda G(x^*)) \big) \ominus_g G(x^*)\\
&\hspace*{-6.5cm}=\lambda G(y) \ominus_g \lambda G(x^*)=\lambda (G(y) \ominus_g  G(x^*))
\end{align*}
So,
$$ - \frac{F(x^* +\lambda (y-x^*)) \ominus_g F(x^*)}{\lambda} \leq G(y) \ominus_g G(x^*) $$
Taking the limit as $ \lambda \to 0^+$ in the last inequality,
$$ -F'(x^*;y-x^*) \leq G(y) \ominus_g G(x^*). $$
By Theorem \ref{t},
$$ -\nabla F(x^*) (y-x^*) \leq G(y) \ominus_g G(x^*). $$
i.e., $ -\nabla F(x^*) \in \partial G(x^*).$
(II) Let $F$ be convex. According to Theorem \ref{x}, $x^*$ is a local minimum solution of $(P)$ if and only if $x^*$ is a global minimum. If $x^*$ is an optimal solution of $(P),$ then as we proved in I, \eqref{AN} is satisfied. Now, suppose that \eqref{AN} holds, then we show that $x^*$ is a global optimal solution of $(P).$\\
By using \eqref{AN}, we have
\begin{equation}\label{31}
 G(y)\ominus_g G(x^*) \geq -\nabla F(x^*)(y-x^*).
\end{equation}
Also using Theorem \ref{tt}, we get
$$ F(y) \ominus_g F(x^*) \geq F'(x^*;y-x^*).$$
Applying Theorem \ref{t}, we obtain
\begin{equation}\label{32}
F(y) \ominus_g F(x^*) \geq \nabla F(x^*)(y-x^*),
\end{equation}
for any $ y \in M.$ Adding \eqref{31} and \eqref{32}, we have
$$ (G(y) \ominus_g G(x^*)) + (F(y) \ominus_g F(x^*)) \geq \hat{0}.$$
As a result,
$$(G(y) + F(y)) \ominus_g (G(x^*) + F(x^*)) \geq \hat{0}.$$
Therefore,
$$ G(y) + F(y) \geq G(x^*) + F(x^*),$$
for any $ y \in M.$ That is, $x^*$ is an optimal solution of $(P).$
\end{proof}
\section{\bf Conclusion}
In this paper, we have extended the concept of directional derivative for fuzzy $n$-cell number-valued functions by $g$-difference \cite{Bede,Gomes1}, and we have given important theorems related to it. We have proposed a definition of the subdifferential of fuzzy $n$-cell number valued functions based on the concept of $g$-difference and partial order defined on the set of n-dimensional fuzzy numbers \cite{Wang}. Then, we have tried to prove theorems related to the subdifferential with this new definition. Our motivation for introducing this new definition is that if we compute the subdifferential, then we can optimize any convex function. Therefore, we proved some optimization theorems using subdifferential and discussed the dual Lagrange problems and Karush–Kuhn–Tucker optimality conditions. Moreover, some necessary and sufficient conditions for the optimal solution of a composite function by subgradient are studied.

\end{document}